\documentclass[notitlepage,leqno,10pt]{amsart}

\usepackage{color}
\usepackage{amsfonts}
\usepackage{amsmath}
\usepackage{amssymb}
\usepackage{supertabular}
\usepackage[latin1]{inputenc}

\newcommand{\cb}{\color{blue}}

\newcommand{\R}{\mathbb{R}}

\newcommand{\la}{\lambda}

\renewcommand{\proof}[1]{\par\smallskip\noindent{\bf Proof#1.}}
\renewcommand{\qed}{\penalty 500\hfill$\square$\par\medskip}

\newcommand{\td}{\tilde}

\newcommand{\ov}{\overline}

\newcommand{\pa}{\partial}

\newcommand{\avg}{\raise 0pt\hbox{$-$}\hskip -10.7pt\int}
\newcommand{\pd}[2]{{\frac{\pa #1}{\pa #2}}}

\newcommand{\dist}{\mathop{\rm dist}}

\newcommand{\eps}{\varepsilon}
\def\epsilon{\eps}

\newcommand{\restrict}[1]{|_{\raise-2pt\hbox{$\scriptstyle #1$}}}
\newcommand{\eq}[1]{\begin{equation}\label{#1}}
\newcommand{\eeq}{\end{equation}}

\newtheorem{lem} {Lemma}[section]
\newtheorem{pro} [lem]{Proposition}
\newtheorem{thm} [lem]{Theorem}
\newtheorem{rem} [lem]{Remark}

\begin{document}

\title[Stable solutions for $\Delta^2u = \lambda e^u$]
{Stable solutions for the bilaplacian with exponential nonlinearity.}

\author{Juan D\'avila}
\address{Departamento de Ingenier\'{\i}a Matem\'atica and CMM,
Universidad de Chile, Casilla 170 Correo 3, Santiago, Chile.}
\email{jdavila@dim.uchile.cl}

\author{Louis Dupaigne}
\address{Laboratoire Amienois de Mathematique
Fondamentale e Appliquee \\ Faculte de Mathematique et
d'Informatique \\ Amiens, France.}
\email{louis.dupaigne@u-picardie.fr}

\author{Ignacio Guerra}
\address{Departamento de Matematica y C.C., Facultad de Ciencia \\
Universidad de Santiago de Chile, Casilla 307, Correo 2, Santiago,
Chile.} \email{iguerra@usach.cl}

\author{Marcelo Montenegro}
\address{Universidade Estadual de
Campinas, IMECC, Departamento de Matem\'atica, \\
Caixa Postal 6065, CEP 13083-970, Campinas, SP, Brasil.}
\email{msm@ime.unicamp.br}

\date{}

\subjclass{Primary 35J65, Secondary 35J40 } \keywords{Biharmonic,
singular solutions, stability.}

\begin{abstract}

Let $\lambda^*>0$ denote the largest possible value of $\lambda$
such that
\begin{align*}
\left\{
\begin{aligned}
\Delta^2 u & = \la e^u && \text{in $B $ }
\\
u &= \pd{u}{n} = 0 && \text{on $ \pa B $ }
\end{aligned}
\right.
\end{align*}
has a solution, where $B$ is the unit ball in $\R^N$ and $n$ is
the exterior unit normal vector. We show that for
$\lambda=\lambda^*$ this problem possesses a unique {\em weak}
solution $u^*$. We prove that $u^*$ is smooth if $N\le 12$ and
singular when $N\ge 13$, in which case $ u^*(r) = - 4 \log r +
\log ( 8(N-2)(N-4) / \lambda^*) + o(1)$ as $r\to 0$. We also
consider the problem with general constant Dirichlet boundary
conditions.
\end{abstract}

\setcounter{page}{1}

\maketitle

%\noindent {\sc AMS Mathematics Subject Classification 2000: 35J65
%35J40 }:

%%%%%%%%%%%%%%%%%%%%%%%%%%%%%%%%%%%%%%%%%%%%%%%%%%%%%%%%%%%%%%%%%%%%%%

\section{Introduction}

We study the fourth order problem
\begin{align}\label{main2}
\left\{
\begin{aligned}
\Delta^2 u & = \la e^u && \text{in $B $ }
\\
u &= a && \text{on $ \pa B $ }
\\
\pd{u}{n} &= b && \text{on $ \pa B $}
\end{aligned}
\right.
\end{align}
where $a$, $b\in \R$, $ B$ is the unit ball in $\R^N$, $N \ge 1$,
$n$ is the exterior unit normal vector and $\la \ge 0$ is a
parameter.

Recently higher order equations have attracted the interest of
many researchers. In particular fourth order equations with an
exponential non-linearity have been studied in 4 dimensions, in a
setting analogous to Liouville's equation, in
\cite{baraket-dammak-pacard-ouni,djadli-malchiodi,wei} and in
higher dimensions by \cite{agg,aggm, bg, bgm,Ferrero-Grunau}.

We shall pay special attention to \eqref{main2} in the case
$a=b=0$, as it is the natural fourth order analogue of the
classical Gelfand problem
\begin{align}\label{250}
\left\{
\begin{aligned}
-\Delta u &= \lambda e^u && \hbox{in $\Omega$}
\\
u&=0 && \hbox{on $\partial \Omega$}
\end{aligned}
\right.
\end{align}
($\Omega$ is a smooth bounded domain in $\R^N$) for which a vast
literature exists
\cite{bcmr,bv,ck,CR,jl,m,mignot-puel,mignot-puel2}.

From the technical point of view, one of the basic tools in the
analysis of \eqref{250} is the maximum principle. As pointed out
in \cite{aggm}, in general domains the maximum principle for
$\Delta^2$ with Dirichlet boundary condition is not valid anymore.
One of the reasons to study \eqref{main2} in a ball is that a
maximum principle holds in this situation, see \cite{boggio}. In
this simpler setting, though there are some similarities between
the two problems, several tools that are well suited for
\eqref{250} no longer seem to work for \eqref{main2}.

As a start, let us introduce the class of weak solutions we shall
be working with~: we say that $u \in H^2(B)$ is a weak solution to
\eqref{main2} if $e^u \in L^1(B)$, $u = a$ on $\partial B$,
$\pd{u}{n}=b$ on $\partial B$ and
\begin{align*}
\int_B \Delta u \Delta \varphi = \la \int_B e^u \varphi, \quad
\hbox{for all } \varphi \in C_0^\infty(B).
\end{align*}

The following basic result is a straightforward adaptation of
Theorem 3 in \cite{aggm}.
\begin{thm}{\rm(\cite{aggm})}\label{t0}
There exists $\la^*$ such that if $0 \le \la < \la^*$ then
\eqref{main2} has a minimal smooth solution $u_\la $ and if $\la >
\la^*$ then \eqref{main2} has no weak solution.

The limit $ u^* = \lim_{\la \nearrow \la^*} u_\la$ exists
pointwise, belongs to $H^2(B)$ and is a weak solution to
\eqref{main2}. It is called the extremal solution.

The functions $u_\la$, $0 \le \la < \la^*$ and $u^*$ are radially
symmetric and radially decreasing.
\end{thm}

The branch of minimal solutions of {\eqref{main2}} has an
important property, namely $u_\lambda$ is stable in the sense that
\begin{align}\label{stability}
\int_B (\Delta \varphi)^2 \ge \lambda \int_B e^{u_\lambda}
\varphi^2, \quad \forall \varphi \in C_0^\infty(B),
\end{align}
see \cite[Proposition 37]{aggm}.

The authors in \cite{aggm} pose several questions, some of which
we address in this work.

\noindent First we show that the extremal solution $u^*$ is the
unique solution to \eqref{main2} in the class of weak solutions.
Actually the statement is stronger, asserting that for
$\lambda=\lambda^*$ there are no strict super-solutions.
\begin{thm}\label{t2} If
\begin{align}\label{t2-1}
\hbox{$v \in H^2(B)$, $e^v \in L^1(B)$, $v|_{\partial B}=a$,
$\frac{\partial v}{\partial n}|_{\partial B} \le b$}
\end{align}
and
\begin{align}\label{t2-2}
\int_B \Delta v \Delta \varphi \ge \lambda^* \int_B e^v \varphi
\quad \forall\;\varphi \in C_0^\infty(B), \,\varphi\ge0,
\end{align}
then $v=u^*$. In particular for $\lambda = \lambda^*$ problem
\eqref{main2} has a unique weak solution.
\end{thm}
This result is analogous to work of Martel \cite{m} for more
general versions of \eqref{250} where the exponential function is
replaced by a positive, increasing, convex and superlinear
function.

Next, we discuss the regularity of the extremal solution $u^*$. In
dimensions $N=5,\ldots,16$ the authors of \cite{aggm} find, with a
computer assisted proof, a radial singular solution $U_\sigma$ to
\eqref{main2} with $a=b=0$ associated to a parameter $\la_\sigma>
8 (N-2) (N-4)$. They show that $\lambda_\sigma < \lambda^*$ if $N
\le 10$ and claim to have numerical evidence that this holds for
$N\le 12$. They leave open the question of whether $u^*$ is
singular in dimension $N\le 12$. We prove
\begin{thm}\label{t1}
If $N \le 12 $ then the extremal solution $u^*$ of \eqref{main2}
is smooth.
\end{thm}
The method introduced in \cite{CR,mignot-puel} to prove the
boundedness of $u^*$ in low dimensions for \eqref{250} seems not
useful for \eqref{main2}, thus requiring a new strategy. A first
indication that the borderline dimension for the boundedness of
$u^*$ is 12 is Rellich's inequality \cite{rellich}, which states
that if $N\ge 5$ then
\begin{align}\label{hardy-rellich}
\int_{\R^N} ( \Delta \varphi)^2 \ge \frac{N^2 (N-4)^2}{16}
\int_{\R^N} \frac{\varphi^2}{|x|^4} \qquad \forall \varphi \in
C_0^\infty(\R^N),
\end{align}
where the constant $N^2 (N-4)^2 / 16$ is known to be optimal. The
proof of Theorem~\ref{t1} is based on the observation that if
$u^*$ is singular then $\lambda^* e^{u^*} \sim 8(N-2)(N-4)
|x|^{-4}$ near the origin. But $8(N-2)(N-4) > N^2 (N-4)^2 /16$ if
$N\le 12$ which would contradict the stability condition
\eqref{stability}.

In view of Theorem~\ref{t1}, it is natural to ask whether $u^*$ is
singular in dimension $N \ge 13$. If $a=b=0$, we prove
\begin{thm}\label{t3}
Let $N \ge 13$ and $a=b=0$. Then the extremal solution $u^*$ to
\eqref{main2} is unbounded.
\end{thm}

For general boundary values, it seems more difficult to determine
the dimensions for which the extremal solution is singular. We
observe first that given any $a,b\in\R$, $u^*$ is the extremal
solution of \eqref{main2} if and only if $u^*-a$ is the extremal
solution of the same equation with boundary condition $u=0$ on
$\partial B$. In particular, if $\lambda^*(a,b)$ denotes the
extremal parameter for problem \eqref{main2}, one has that
$\lambda^*(a,b)=e^{-a}\la^*(0,b)$. So the value of $a$ is
irrelevant. But one may ask if Theorem~\ref{t3} still holds for
any $N\ge 13$ and any $b\in\R$. The situation turns out to be
somewhat more complicated~:

\begin{pro}\label{p1.7}

\

\begin{itemize}
\item[a)] Fix $N\ge 13$ and take any $a\in\R$. Assume $b\ge -4$. There exists a
critical parameter $b^{max}>0$, depending only on $N$, such that the
extremal solution $u^*$ is singular if and only if $b\le b^{max}$.

\item[b)] Fix $b \ge -4$ and take any $a\in\R$. There exists a critical dimension
$N^{min}\ge13$, depending only on $b$, such that the extremal
solution $u^*$ to \eqref{main2} is singular if $N\ge N^{min}$.
\end{itemize}
\end{pro}

\begin{rem}

\

\begin{itemize}
\item We have not investigated the case $b<-4$. \item If follows
from item {\rm a)} that for $b\in[-4,0]$, the extremal solution is
singular if and only if $N\ge13$. \item It also follows from item
{\rm a)} that there exist values of $b$ for which $N^{min}>13$. We
do not know whether $u^*$ remains bounded for $13\le N<N^{min}$.
\end{itemize}
\end{rem}

Our proof of Theorem~\ref{t3} is related to an idea that Brezis
and V\'azquez applied for the Gelfand problem and is based on a
characterization of singular {\em energy} solutions through
linearized stability (see Theorem~3.1 in \cite{bv}). In our
context we show
\begin{pro}\label{240}
Assume that $u \in H^2(B)$ is an unbounded weak solution of
{\eqref{main2}} satisfying the stability condition
\begin{align}\label{260}
\la \int_{B} e^{u} \varphi^2 \le \int_{B} ( \Delta \varphi)^2,
\quad \forall \varphi \in C_0^\infty(B).
\end{align}
Then $\lambda = \lambda^*$ and $u=u^*$.
\end{pro}

We do not use Proposition \ref{240} directly but some variants of
it -- see Lemma~\ref{comp-non-linear} and
Remark~\ref{comp-non-linear-2} in Section~\ref{270} -- because we
do not have at our disposal an explicit solution to the equation
\eqref{main2}. Instead, we show that it is enough to find a
sufficiently good approximation to $u^*$. When $N \ge 32$ we are
able to construct such an approximation by hand. However, for
$13\le N \le 31$ we resort to a computer assisted generation and
verification.

Only in very few situations one may take advantage of
Proposition~\ref{240} directly. For instance for problem
\eqref{main2} with $a=0$ and $b=-4$ we have an explicit solution
\begin{align*}
\bar u(x) = - 4 \log |x|
\end{align*}
associated to $\bar\lambda = 8(N-2)(N-4)$. Thanks to Rellich's
inequality \eqref{hardy-rellich} the solution $\bar u$ satisfies
condition \eqref{260} when $N \ge 13$. Therefore, by
Theorem~\ref{t1} and a direct application of Proposition~\ref{240}
we obtain Theorem~\ref{t3} in the case $b=-4$.

In \cite{aggm} the authors say that a radial weak solution $u$ to
\eqref{main2} is weakly singular if
\begin{align*}
\lim_{r\to 0} r u'(r) \quad { exists.}
\end{align*}
For example, the singular solutions $U_\sigma$ of \cite{aggm}
verify this condition.

As a corollary of Theorem~\ref{t2} we show
\begin{pro}\label{weakly singular}
{The extremal solution $u^*$ to \eqref{main2}} with $b \ge -4$ is
always weakly singular.
\end{pro}
A weakly singular solution either is smooth or exhibits a log-type
singularity at the origin. More precisely, if $u$ is a non-smooth
weakly singular solution of \eqref{main2} with parameter $\lambda$
then (see \cite{aggm})
\begin{align*}
&\lim_{r\to 0} u(r) + 4 \log r = \log\frac{8(N-2)(N-4)}{\lambda},
\\
&\lim_{r\to 0} r u'(r) = -4.
\end{align*}

In Section~\ref{270} we describe the comparison principles we use
later on. Section~\ref{s-unique} is devoted to the proof of the
uniqueness of $u^*$ and Propositions~\ref{240} and \ref{weakly
singular}. We prove Theorem~\ref{t1}, the boundedness of $u^*$ in
low dimensions, in Section~\ref{s-bounded-u}. The argument for
Theorem~\ref{t3} is contained in Section~\ref{s-unbounded-1} for
the case $N\ge 32$ and Section~\ref{s-unbounded-2} for $13\le N
\le 31$. In Section~\ref{s-last} we give the proof of
Proposition~\ref{p1.7}.

\bigskip\noindent
{\bf Notation.}
\begin{itemize}
\item $ B_R $: ball of radius $R$ in $\R^N$ centered at the
origin. $B = B_1$. \item $n$: exterior unit normal vector to $B_R$
%\item
%$ a^+ = \max(a,0 )$ for $a\in\R$.
\item All inequalities or equalities for functions in $L^p$ spaces
are understood to be a.e.
\end{itemize}

%%%%%%%%%%%%%%%%%%%%%%%%%%%%%%%%%%%%%%%%%%%%%%%%%%%%%%%%%%%%%%%%%%%%%%

\section{Comparison principles}\label{270}

\begin{lem}(Boggio's principle, \cite{boggio})
If $u \in C^4(\ov B_R)$ satisfies
\begin{equation*}
\left\{
\begin{aligned}
\Delta^2 u & \ge 0 && \text{in $B_R $ }
\\
u &= \pd{u}{n} = 0 && \text{on $ \pa B_R $ }
\end{aligned}
\right.
\end{equation*}
then $u\ge 0$ in $B_R $.
\end{lem}

\begin{lem}\label{weak-max-prin}
Let $ u \in L^1(B_R)$ and suppose that
\begin{align*}
\int_{B_R} u \Delta^2 \varphi \ge 0
\end{align*}
for all $\varphi \in C^4(\ov B_R)$ such that $\varphi \ge 0$ in $
B_R$, $ \varphi|_{\pa B_R} = 0 = \pd\varphi n |_{\pa B_R}$. Then $
u \ge 0$ in $B_R$. Moreover $u\equiv 0$ or $u>0$ a.e.\@ in $B_R$.
\end{lem}
For a proof see Lemma~17 in \cite{aggm}.

\begin{lem}\label{radial-comp}
If $u \in H^2(B_R)$ is radial, $\Delta^2 u \ge 0$ in $B_R$ in the
weak sense, that is
\begin{align*}
\int_{B_R} \Delta u \Delta \varphi \ge 0 \quad \forall \varphi \in
C_0^\infty(B_R), \, { \varphi \ge0}
\end{align*}
and $u|_{\pa B_R}\ge 0$, $\pd un|_{\pa B_R}\le 0$ then $u \ge 0$
in $B_R$.
\end{lem}
%\bigskip
%\noindent { \color{red} \large \bf I don't know if this holds for
%nonradial functions}
%\bigskip
\proof{} We only deal with the case $R=1$ for simplicity.
%{\color{red}I added this restriction because there was a
%small misprint with the bc's for $u_2$ in the previous version.}
Solve
\begin{align*}
\left\{
\begin{aligned}
\Delta^2 u_1 &= \Delta^2 u && \text{in $B_1 $ }
\\
u_1 &= \pd{u_1}{n} = 0 && \text{on $ \pa B_1 $ }
\end{aligned}
\right.
\end{align*}
in the sense $u_1 \in H_0^2(B_1) $ and $\int_{B_1} \Delta u_1
\Delta \varphi = \int_{B_1} \Delta u \Delta \varphi $ for all $
\varphi \in C_0^\infty(B_1)$. Then $u_1 \ge 0$ in $B_1 $ by
Lemma~\ref{weak-max-prin}.

Let $u_2 = u - u_1 $ so that $\Delta^2 u_2 =0 $ in $B_1 $. Define
$ f = \Delta u_2$. Then $ \Delta f=0 $ in $B_1 $ and since $f$ is
radial we find that $ f$ is constant. It follows that $u_2 = a r^2
+ b$. Using the boundary conditions we deduce $ a+b \ge 0$ and $ a
\le 0$, which imply $ u_2 \ge 0$. \qed

{Similarly we have
\begin{lem}
If $u \in H^2(B_R)$ and $\Delta^2 u \ge 0$ in $B_R$ in the weak
sense, that is
\begin{align*}
\int_{B_R} \Delta u \Delta \varphi \ge 0 \quad \forall \varphi \in
C_0^\infty(B_R), \, \varphi \ge0
\end{align*}
and $u|_{\pa B_R}= 0$, $\pd un|_{\pa B_R}\le 0$ then $u \ge 0$ in
$B_R$.
\end{lem}}

The next lemma is a consequence of a decomposition lemma of Moreau
\cite{moreau}. For a proof see \cite{gazzola-grunau,ggm}.
\begin{lem} \label{moreau}
Let $u\in H_0^2(B_R) $. Then there exist unique $w,v \in
H_0^2(B_R) $ such that $u=w+v$, $w \ge 0 $, $\Delta^2 v \le 0 $ in
$ B_R$ and $\int_{B_R} \Delta w \Delta v =0$.
\end{lem}

{We need the following comparison principle.}
\begin{lem} \label{comp-non-linear}
Let $u_1$, $u_2 \in H^2(B_R)$ with $e^{u_1}$, $e^{u_2} \in
L^1(B_R)$. Assume that
\begin{equation*}
\Delta^2 u_1 \le \la e^{u_1} \quad \hbox{in $B_R $ }
\end{equation*}
in the sense
\begin{align}\label{weak}
\int_{B_R} \Delta u_1 \Delta \varphi \le \la \int_{B_R} e^{u_1}
\varphi \quad \forall \varphi \in C_0^\infty(B_R), \, \varphi \ge
0,
\end{align}
and $\Delta^2 u_2 \ge \la e^{u_2}$ in $B_R $ in the similar weak
sense. Suppose also
\begin{equation*}
u_1|_{\pa B_R} = u_2|_{\pa B_R} \quad\hbox{and}\quad \pd
{u_1}n|_{\pa B_R} = \pd {u_2}n|_{\pa B_R}.
\end{equation*}
Assume furthermore that $u_1$ is stable in the sense that
\begin{equation}\label{stability minimal}
\la \int_{B_R} e^{u_1} \varphi^2 \le \int_{B_R} ( \Delta
\varphi)^2, \quad \forall \varphi \in C_0^\infty(B_R).
\end{equation}
Then
\[
u_1 \le u_2 \quad \hbox{ in $B_R$.}
\]
\end{lem}
\proof{}{Let $u = u_1 - u_2$. By Lemma~\ref{moreau} there exist $w,v
\in H_0^2(B_R)$ such that $u=w+v$, $w\ge 0$ and $\Delta^2 v \le 0$.}
Observe that $v\le 0 $ so $w \ge u_1 - u_2$.

By hypothesis we have for all $\varphi \in C_0^\infty(B_R), \,
\varphi \ge 0$,
\begin{equation*}
\int_{B_R} \Delta (u_1- u_2) \Delta \varphi \le \la \int_{B_R} (
e^{u_1} - e^{u_2} )\varphi \le \la \int_{B_R\cap[u_1\ge u_2]} (
e^{u_1} - e^{u_2} )\varphi
\end{equation*}
and by density this holds also for $w$:
\begin{multline} \label{1}
\int_{B_R} (\Delta w )^2 = \int_{B_R} \Delta (u_1- u_2) \Delta w
\\ \le \la \int_{B_R\cap[u_1\ge u_2]} ( e^{u_1} - e^{u_2} )w = \la
\int_{B_R} ( e^{u_1} - e^{u_2} )w,
\end{multline}
where the first equality holds because $\int_{B_R} \Delta w \Delta
v =0$. By density we deduce from {\cb\eqref{stability minimal}:}
\begin{equation}\label{2}
\la \int_{B_R} e^{u_1} w^2 \le \int_{B_R} (\Delta w)^2.
\end{equation}
%{\color{red} the density result should be standard but I don't
%know a reference}
Combining \eqref{1} and \eqref{2} we obtain
\begin{equation*}
\int_{B_R} e^{u_1} w^2 \le \int_{B_R} ( e^{u_1} - e^{u_2} ) w .
\end{equation*}
Since $u_1 - u_2 \le w $ the previous inequality implies
\begin{equation}\label{3}
0 \le \int_{B_R} ( e^{u_1} - e^{u_2} - e^{u_1} ( u_1 - u_2) ) w.
\end{equation}
But by convexity of the exponential function $ e^{u_1} - e^{u_2} -
e^{u_1} ( u_1 - u_2) \le 0$ and we deduce from \eqref{3} that $ (
e^{u_1} - e^{u_2} - e^{u_1} ( u_1 - u_2) ) w = 0$. Recalling that
$u_1 - u_2 \le w $ we deduce that $ u_1 \le u_2$. \qed

\begin{rem}\label{comp-non-linear-2}
The following variant of Lemma~\ref{comp-non-linear} also holds:

Let $u_1$, $u_2 \in H^2(B_R)$ be {\bf radial} with $e^{u_1}$,
$e^{u_2} \in L^1(B_R)$. Assume $\Delta^2 u_1 \le \la e^{u_1}$ in
$B_R $ in the sense of \eqref{weak} and $\Delta^2 u_2 \ge \la
e^{u_2}$ in $B_R $. Suppose $ u_1|_{\pa B_R} \le u_2|_{\pa B_R} $
and $ \pd {u_1}n|_{\pa B_R} \ge \pd {u_2}n|_{\pa B_R} $ and that
the stability condition \eqref{stability minimal} holds. Then $u_1
\le u_2$ in $B_R$.

\end{rem}
\proof{} We solve for $ \td u \in H_0^2(B_R)$ such that
\begin{align*}
\int_{B_R} \Delta \td u\Delta \varphi = \int_{B_R} \Delta ( u_1 -
u_2) \Delta \varphi \quad \forall \varphi \in C_0^\infty(B_R).
\end{align*}
By Lemma \ref{radial-comp} it follows that $ \td u \ge u_1 - u_2$.
Next we apply the decomposition of Lemma~\ref{moreau} to $\td u$,
that is $\tilde u = w + v $ with {$w,v \in H_0^2(B_R)$}, $w \ge 0
$, $\Delta^2 v \le 0 $ in $ B_R$ and $\int_{B_R} \Delta w \Delta v
=0$. Then the argument follows that of
Lemma~\ref{comp-non-linear}. \qed

Finally, in several places we will need the method of sub and
supersolutions in the context of weak solutions.
\begin{lem}
Let $\lambda> 0$ and assume that there exists $\bar u \in
H^2(B_R)$ such that $e^{\bar u} \in L^1(B_R)$,
\begin{align*}
\int_{B_R} \Delta \bar u \Delta \varphi \ge \la \int_{B_R} e^{\bar
u} \varphi \quad \hbox{for all } \varphi \in C_0^\infty(B_R) ,
\varphi \ge 0
\end{align*}
and
\begin{align*}
{\cb \bar u = a} , \quad \pd{\bar u}{n} \le b \quad \hbox{on
$\partial B_1$}.
\end{align*}
Then there exists a weak solution to \eqref{main2} such that $u
\le \bar u$.
\end{lem}
The proof is similar to that of Lemma 19 in \cite{aggm}.

%%%%%%%%%%%%%%%%%%%%%%%%%%%%%%%%%%%%%%%%%%%%%%%%%%%%%%%%%%%%%%%%%%%%%%
\section{Uniqueness of the extremal solution: proof of
Theorem~\ref{t2}}\label{s-unique}

\proof{ of Theorem~\ref{t2}} Suppose that $v \in H^2(B)$ satisfies
\eqref{t2-1}, \eqref{t2-2} and $v \not \equiv u^*$. Notice that we
do not need $v$ to be radial.

The idea of the proof is as follows :

\medskip \noindent {\bf Step 1.} The function
\begin{align*}
u_0 = \frac{1}{2}(u^* + v)
\end{align*}
is a super-solution to the following problem {
\begin{align} \label{170}
\left\{
\begin{aligned}
\Delta^2 u &= \lambda^* e^u + \mu \eta e^u && \hbox{in $B$}
\\
u &= a && \hbox{on $\partial B$}
\\
\pd un &= b && \hbox{on $\partial B$}
\end{aligned}
\right.
\end{align}}
for some $\mu=\mu_0>0$, where $\eta \in C_0^\infty(B)$, $0\le
\eta\le 1$ is a fixed {radial} cut-off function such that
\begin{align*}
\eta(x)=1 \quad \hbox{for $|x|\le \frac 12$},\qquad \eta(x)=0
\quad \hbox{for $|x| \ge \frac 34$}.
\end{align*}

%\medskip \noindent {\bf Step 2.} For $0<\mu<\mu_0$ equation \eqref{170}
%has a smooth solution.

\medskip \noindent {\bf Step 2.} Using a solution to \eqref{170}
we construct, for some $\lambda>\lambda^*$, a super-solution to
\eqref{main2}. This provides a solution $u_\lambda$ for some
$\lambda>\lambda^*$, which is a contradiction.

\proof{ of Step 1} Observe that given $0<R<1$ we must have for
some $c_0 = c_0(R)>0$
\begin{align}\label{150}
v(x) \ge u^*(x) + c_0 \quad |x|\le R.
\end{align}
To prove this we recall the Green's function for $\Delta^2$ with
Dirichlet boundary conditions {
\begin{align*}
\left\{
\begin{aligned}
\Delta^2_x G (x,y)&= \delta_y && x \in B
\\
G(x,y) &= 0 && x \in \partial B
\\
\pd Gn(x,y) &=0 && x \in \partial B,
\end{aligned}
\right.
\end{align*}
where $\delta_y$ is the Dirac mass at $y \in B$.} Boggio gave an
explicit formula for $G(x,y)$ which was used in
\cite{grunau-sweers} to prove that in dimension $N \ge 5$ (the
case $1\le N\le4$ can be treated similarly)
\begin{align}\label{140}
G(x,y) \sim |x-y|^{4-N} \min\left( 1, \frac{d(x)^2
d(y)^2}{|x-y|^4}\right)
\end{align}
where
\begin{align*}
d(x) = \dist(x,\partial B) = 1 - |x|.
\end{align*}
and $ a\sim b$ means that for some constant $C>0$ we have $C^{-1}
a\le b \le C a$ (uniformly for $x,y\in B$). Formula \eqref{140}
yields
\begin{align} \label{190}
G(x,y) \ge c d(x)^2 d(y)^2
\end{align}
for some $c>0$ and this in turn implies that for smooth functions
$\tilde v$ and $\tilde u$ such that $\tilde v-\tilde u\in
H^2_0(B)$ and $\Delta^2(\tilde v -\tilde u)\ge0$,
\begin{align*}
\tilde v(y) - \tilde u(y) & = \int_{\partial B} \Big( \pd{\Delta_x
G}{n_x}(x,y)
(\tilde v-\tilde u)- \Delta_x G (x,y)\pd{(\tilde v-\tilde u)}n \Big) \, d x \\
& \qquad + \int_B G(x,y) \Delta^2 (\tilde v-\tilde u) \, d x\\
& \ge c d(y)^2 \int_B (\Delta^2{\tilde v} - \Delta^2{\tilde u}
)d(x)^2 \, d x .
\end{align*}
Using a standard approximation procedure, we conclude that
\begin{align*}
v(y) - u^*(y) \ge c d(y)^2 \lambda^* \int_B (e^{v} - e^{u^*}
)d(x)^2 \, d x .
\end{align*}
Since $v \ge u^*$, $v\not\equiv u^*$ we deduce \eqref{150}.

Let $u_0 = (u^* + v)/2$. Then by Taylor's theorem
\begin{align}\label{120}
e^{v} = e^{u_0} + (v - u_0 ) e^{u_0} + \frac12 (v - {u_0})^2
e^{u_0} + \frac16 (v - {u_0})^3 e^{u_0} + \frac1{24} (v - {u_0})^4
e^{\xi_2}
\end{align}
for some ${u_0} \le \xi_2 \le v$ and
\begin{align}\label{130}
e^{u^*} = e^{u_0} + (u^* - {u_0}) e^{u_0} + \frac12 (u^* -
{u_0})^2 e^{u_0} + \frac16 (u^* - {u_0})^3 e^{u_0} + \frac1{24}
(u^* - {u_0})^4 e^{\xi_1}
\end{align}
for some $u^* \le \xi_1 \le {u_0}$. Adding \eqref{120} and
\eqref{130} yields
\begin{align}\label{160}
\frac12(e^{v} + e^{u^*}) \ge e^{u_0} + \frac{1}{8}(v-u^*)^2
e^{u_0}.
\end{align}
From \eqref{150} with $R=3/4$ and \eqref{160} we see that ${u_0} =
(u^* + v)/2$ is a super-solution of \eqref{170} with $\mu_0 :=
c_0/8$.

\proof{ of Step 2} Let us show now how to obtain a weak
super-solution of \eqref{main2} for some $\lambda>\lambda^*$.
Given $\mu>0$, let $u$ denote the minimal solution to \eqref{170}.
Define $\varphi_1$ as the solution to {
\begin{align*}
\left\{
\begin{aligned}
\Delta^2 \varphi_1 &= \mu \eta e^u && \hbox{in $B$}
\\
\varphi_1 &= 0 && \hbox{on $\partial B$}
\\
\pd{\varphi_1}n &= 0 && \hbox{on $\partial B$},
\end{aligned}
\right.
\end{align*}
and $\varphi_2$ be the solution of
\begin{align*}
\left\{
\begin{aligned}
\Delta^2 \varphi_2 &= 0 && \hbox{in $B$}
\\
\varphi_2 &= a && \hbox{on $\partial B$}
\\
\pd{\varphi_2} n &= b && \hbox{on $\partial B$}.
\end{aligned}
\right.
\end{align*}}
If $N\ge5$ (the case $1\le N\le 4$ can be treated similarly),
relation \eqref{190} yields
\begin{align}\label{200}
\varphi_1(x) \ge c_1 d(x)^2 \quad \hbox{for all $x\in B$},
\end{align}
for some $c_1>0$. {But $u$ is a radial solution of \eqref{170} and
therefore it is smooth in $B\setminus B_{1/4}$. Thus
\begin{align}\label{210}
u(x) \le M \varphi_1 + \varphi_2 \quad\hbox{for all $x\in B_{1/2}
$},
\end{align}
for some $M>0$. Therefore, from \eqref{200} and \eqref{210}, for
$\lambda>\lambda^*$ with $\lambda-\lambda^*$ sufficiently small we
have
\begin{align*}
({\textstyle \frac\lambda{\lambda^*}} -1 ) u \le \varphi_1 +
({\textstyle \frac\lambda{\lambda^*}} -1 ) \varphi_2 \quad
\hbox{in $B$.}
\end{align*}
Let {$w = \frac\lambda{\lambda^*} u - \varphi_1 - ({\textstyle
\frac\lambda{\lambda^*}} -1 ) \varphi_2$}. The inequality just
stated guarantees that $w \le u$. Moreover
\begin{align*}
\Delta^2 w = \lambda e^u + \frac{\lambda \mu}{\lambda^*} \eta e^u
- \mu \eta e^u \ge \lambda e^u \ge \lambda e^w \quad \hbox{ in } B
\end{align*}
and
\begin{align*}
w=a \quad \pd{w}{n}=b \quad \hbox{on } \partial B.
\end{align*}
Therefore $w$ is a super-solution to \eqref{main2} for $\lambda$.
By the method of sub and super-solutions a solution to
\eqref{main2} exists for some $\lambda>\lambda^*$, which is a
contradiction. \qed

\proof{ of Proposition~\ref{240}} Let $u \in H^2(B)$, $\lambda>0$
be a weak unbounded solution of \eqref{main2}. If
$\lambda<\lambda^*$ from Lemma~\ref{comp-non-linear} we find that
$u \le u_\lambda$ where $u_\lambda$ is the minimal solution. This
is impossible because $u_\lambda$ is smooth and $u$ unbounded. If
$\lambda = \lambda^*$ then necessarily $ u=u^*$ by
Theorem~\ref{t2}.\qed

\proof{ of Proposition~\ref{weakly singular}} Let $u$ denote the
extremal solution of \eqref{main2} with $b\ge -4$. If $u$ is
smooth, then the result is trivial. So we restrict to the case
where $u$ is singular. By Theorem~\ref{t1} we have in particular
that $N\ge 13$. We may also assume that $a=0$. If $b=-4$ by
Theorem~\ref{t2} we know that if $N\ge 13$ then $u = - 4 \log |x|$
so that the desired conclusion holds. Henceforth we assume $b>-4$
in this section.

For $\rho>0$ define
\begin{align*}
u_\rho(r) = u(\rho r) + 4 \log \rho,
\end{align*}
so that
\begin{align*}
\Delta^2 u_\rho = \lambda^* e^{u_\rho} \qquad \hbox{ in }
B_{1/\rho}.
\end{align*}
%Therefore $v_\rho = \frac{d u_\rho}{d \rho}$ solves
%\begin{align*}
%\Delta^2 v_\rho = \lambda^* e^{u^*} v_\rho \qquad \hbox{ in } B_{1/\rho}.
%\end{align*}
Then
\begin{align*}
\frac{d u_\rho}{d \rho}\Big|_{\rho=1,r=1} = u'(1) + 4 > 0.
\end{align*}
Hence, there is $\delta>0$ such that
\begin{align*}
u_\rho(r) < u(r) \quad \hbox{for all } 1-\delta<r\le 1, \,
1-\delta<\rho\le 1.
\end{align*}
This implies
\begin{align} \label{310}
u_\rho(r) < u(r) \quad \hbox{for all } 0<r\le 1, \,
1-\delta<\rho\le 1.
\end{align}
Otherwise set
\begin{align*}
r_0 = \sup \, \{ \, 0<r<1\, | \, u_\rho(r) \ge u(r)\, \}.
\end{align*}
This definition yields
\begin{align}\label{300}
u_\rho(r_0) = u(r_0) \quad \hbox{ and } \quad u_\rho'(r_0) \le
u'(r_0).
\end{align}
Write $\alpha = u(r_0)$, $\beta = u'(r_0)$. Then $u$ satisfies
{
\begin{align} \label{a3} \left\{
\begin{aligned}
\Delta ^2 u & = \lambda e^u && \hbox{ on $B_{r_0}$}
\\
u(r_0) & = \alpha
\\
u'(r_0) &= \beta .
\end{aligned}
\right.
\end{align}
Observe that $u$ is an unbounded $H^2(B_{r_0})$ solution to
\eqref{a3}, which is also stable. Thus Proposition~\ref{240} shows
that $u$ is the extremal solution to this problem. On the other
hand $u_\rho$ is a supersolution to \eqref{a3}, since
$u_\rho'(r_0) \le \beta$ by \eqref{300}.  We may now use
Theorem~\ref{t2} and we deduce that}
\begin{align*}
u(r) = u_\rho(r) \quad \hbox{for all } 0<r \le r_0,
\end{align*}
which in turn implies by standard ODE theory that
\begin{align*}
u(r) = u_\rho(r) \quad \hbox{for all } 0<r \le 1,
\end{align*}
a contradiction {\cb with \eqref{310}}. This proves estimate
\eqref{310}.

From \eqref{310} we see that
\begin{align}\label{315}
\frac{d u_\rho}{d \rho}\Big|_{\rho=1}(r) \ge 0 \quad \hbox{for all
} 0<r\le 1.
\end{align}
But
\begin{align*}
\frac{d u_\rho}{d \rho}\Big|_{\rho=1}(r) = u'( r) r + 4 \quad
\hbox{for all } 0<r\le 1
\end{align*}
and this together with \eqref{315} implies
\begin{align} \label{330}
\frac{d u_\rho}{d \rho}(r) = \frac{1}{\rho}( u'(\rho r) \rho r + 4
) \ge 0 \quad \hbox{for all } 0<r\le \frac{1}{\rho} , \, 0<\rho\le
1.
\end{align}
which means that $u_\rho(r)$ is non-decreasing in $\rho$. We wish
to show that $\lim_{\rho \to 0} u_\rho(r)$ exists for all $0<r\le
1$. For this we shall show
\begin{align}\label{320} u_\rho(r) \ge -
4\log(r) + \log \left(\frac{8(N-2)(N-4)}{\lambda^*}\right) \quad
\hbox{for all }0<r\le \frac1\rho , \, 0<\rho\le 1.
\end{align}
Set
\begin{align*}
u_0(r) = - 4\log(r) + \log
\left(\frac{8(N-2)(N-4)}{\lambda^*}\right).
\end{align*}
and suppose that \eqref{320} is not true for some $0<\rho<1$. Let
\begin{align*}
r_1 = \sup \, \{ \, 0<r<1/\rho \, | \, u_\rho(r) < u_0(r)\, \}.
\end{align*}
Observe that
\begin{align}
\label{aaa10} \lambda^*>8(N-2)(N-4).
\end{align}
Otherwise $w=-4\ln r$ would be a strict supersolution of the
equation satisfied by $u$, which is not possible by Theorem
\ref{t2}. In particular, $r_1<1/\rho$ and
\begin{align*}
u_\rho(r_1) = u_0(r_1) \quad \hbox{ and } \quad u_\rho'(r_1) \ge
u_0'(r_1).
\end{align*}
It follows that $u_0$ is a supersolution of
\begin{align}\label{eqbr1}
\left\{
\begin{aligned}
\Delta^2 u & = \la^* e^u && \text{in $B_{r_1} $ }
\\
u &= A && \text{on $ \pa B_{r_1} $ }
\\
\pd{u}{n} &= B && \text{on $ \pa B_{r_1} $},
\end{aligned}
\right.
\end{align}
with $A=u_\rho(r_1)$} and $B=u_\rho'(r_1)$. Since $u_\rho$ is a
singular stable solution of \eqref{eqbr1}, it is the extremal
solution of the problem by Proposition \ref{240}. By Theorem
\ref{t2}, there is no strict supersolution of \eqref{eqbr1} and we
conclude that $u_\rho\equiv u_0$ first for $0<r<r_1$ and then for
$0<r\le 1/\rho$. This is impossible for $\rho>0$ because
{\cb$u_\rho(1/ \rho) = 4 \log \rho$ and $u_0(1/\rho) < 4 \log \rho
+ \log(\frac{8(N-2)(N-4)}{\lambda^*}) < u_\rho(1/ \rho)$ by
\eqref{aaa10}}. This proves \eqref{320}.

By \eqref{330} and \eqref{320} we see that
\begin{align*}
v (r) = \lim_{\rho\to 0} u_\rho(r) \quad \hbox{ exists for all }
0< r< + \infty,
\end{align*}
where the convergence is uniform (even in $C^k$ for any $k$) on
compact sets of $\R^N \setminus \{0\}$. Moreover $v$ satisfies
\begin{align}\label{340}
\Delta^2 v = \lambda^* e^v \quad \hbox{ in } \R^N\setminus\{0\}.
\end{align}
Then for any $r>0$
\begin{align*}
v(r) = \lim_{\rho\to 0} u_\rho(r) = \lim_{\rho\to 0} u(\rho r) + 4
\log(\rho r) - 4 \log(r) = v(1) - 4 \log (r).
\end{align*}
Hence, using equation \eqref{340} we obtain
\begin{align*}
v(r) = -4 \log r + \log \left(\frac{8(N-2)(N-4)}{\lambda^*}\right)
= u_0(r).
\end{align*}
But then
\begin{align*}
u_\rho'(r) = u'(\rho r) \rho \to - 4, \quad \hbox{ as } \rho \to
0.
\end{align*}
and therefore, with $r=1$
\begin{align}\label{1400}
\rho u'(\rho) \to -4 \quad \hbox{ as } \rho \to 0.
\end{align}
\qed

%%%%%%%%%%%%%%%%%%%%%%%%%%%%%%%%%%%%%%%%%%%%%%%%%%%%%%%%%%%%%%%%%%%%%%

\section{Proof of Theorem~\ref{t1}}\label{s-bounded-u}

We will show first
\begin{lem}
Suppose that the extremal solution $u^*$ to \eqref{main2} is
singular. Then for any $\sigma>0$ there exists $0<R<1$ such that
\begin{align}\label{cotainf}
u^*(x) \ge (1-\sigma) \log\left(\frac{1}{|x|^4}\right) , \quad
\forall\; |x|<R.
\end{align}
\end{lem}
\proof{} Assume by contradiction that \eqref{cotainf} is false.
Then there exists $\sigma>0$ and a sequence $x_k \in B$ with $x_k
\to 0$ such that
\begin{align}\label{absurd}
u^{*}(x_k) < (1-\sigma) \log\left(\frac{1}{|x_k|^4}\right).
\end{align}
Let $s_k = |x_k|$ and choose $0<\lambda_k<\lambda^*$ such that

\begin{align}\label{lambdak}
\max_{\ov B} u_{\lambda_k} = u_{\lambda_k}(0) =
\log\left(\frac1{s_k^4}\right).
\end{align}
Note that $\lambda_k \to \lambda^*$, otherwise $u_{\lambda_k}$
would remain bounded. Let
\begin{align*}
v_k(x) = \frac{u_{\lambda_k}(s_k x)}{\log(\frac1{s_k^4})} \qquad x
\in B_k \equiv\frac1{s_k} B.
\end{align*}
Then $0 \le v_k \le 1$, $v_k(0)=1$,
\begin{align*}
\Delta^2 v_k(x) & = \lambda_k \frac{s_k^4}{\log(\frac{1}{s_k^4}) }
e^{u_{\lambda_k}(s_k x)}
\\
& \le \frac{\lambda_k}{\log(\frac{1}{s_k^4})} \to 0 \quad \hbox{in
$B_k$}
\end{align*}
by \eqref{lambdak}. By elliptic regularity $v_k \to v$ uniformly
on compact sets of $\R^N$ to a function $v$ satisfying $0 \le v
\le 1$, $v(0)=1$, $\Delta^2 v=0$ in $\R^N$.By Liouville's theorem
for biharmonic functions \cite{huilgol} we conclude that $v$ is
constant and therefore $v\equiv 1$.

Since $|x_k| = s_k$ we deduce that
\begin{align*}
\frac{u_{\lambda_k}(x_k)}{\log(\frac1{s_k^4})} \to 1,
\end{align*}
which contradicts \eqref{absurd}. \qed

\proof{ of Theorem~\ref{t1}} We write for simplicity $u = u^*$,
$\lambda = \lambda^*$. Assume by contradiction that $u^*$ is
unbounded and $5 \le N\le12$. If $N \le 4$ the problem is
subcritical, and the boundedness of $u^*$ can be proved by other
means : no singular solutions exist for positive $\lambda$ (see
\cite{aggm})-though in dimension $N=4$ they can blow up as
$\lambda \to 0$, see \cite{wei}.

{For $\eps>0$ let $\psi = |x|^{\frac{4-N}{2}+\eps}$ and let $\eta
\in C_0^\infty(\R^N)$ with $\eta \equiv 1$ in $B_{1/2}$ and
$supp(\eta) \subseteq B$. Observe that
\begin{align*}
(\Delta \psi)^2 = (H_N + O(\eps) ) |x|^{-N + 2 \eps} , \quad
\hbox{where } H_N = \frac{N^2 (N-4)^2}{16} .
\end{align*}
Using a standard approximation argument as in the proof of
Lemma~\ref{comp-non-linear}, we can use $\psi \eta $ as a test
function in \eqref{stability minimal} and we obtain
\begin{align*}
\int_{B} (\Delta \psi)^2 + O(1) \ge \lambda \int_{B} e^u \psi^2,
\end{align*}
since the contribution of the integrals outside a fixed ball
around the origin remains bounded as $\eps \to 0$ (here $O(1)$
denotes a bounded function  as $\eps\to 0$).

This implies
\begin{align}\label{upper estimate}
\lambda \int_{B} e^u |x|^{4-N+2 \eps} \le (H_N + O(\eps) ) \int_B
|x|^{-N+2\eps} =\omega_N \frac{H_N}{2\eps} + O(1)
\end{align}
where $\omega_N$ is the surface area of the unit $N-1$ dimensional
sphere $S^{N-1}$. In particular $\int_{B} e^u |x|^{4-N+2 \eps} < +
\infty$.

}

For $\eps>0$ we define $\varphi=|x|^{4-N+ 2 \eps}$. Note that away
from the origin
\begin{align}\label{bilaplace varphi}
\Delta^2 \varphi = \eps k_N |x|^{-N+2\eps}, \quad \hbox{ where } k_N
= 4 ( N-2)(N-4) + O(\eps).
\end{align}
Let $\varphi_j$ solve
\begin{align}
\label{bilaplace varphi j} \left\{
\begin{aligned}
\Delta^2 \varphi_j &= \eps k_N \min(|x|^{-N+2\eps},j) && \hbox{in
$B$}
\\
\varphi_j &= \pd{\varphi_j}{n} = 0 && \hbox{on $\partial B$.}
\end{aligned}
\right.
\end{align}
Then $\varphi_j \uparrow \varphi$ as $j\to + \infty$. Using
\eqref{upper estimate} and \eqref{bilaplace varphi j}
\begin{align} \nonumber
\eps k_N \int_B u \; \min(|x|^{-N+2\eps} , j) &= \int_B u \Delta^2
\varphi_j = \lambda \int_{B} e^u \varphi_j
\\
\nonumber &\le \lambda \int_{B} e^u \varphi
\\
\nonumber & \le \omega_N \frac{H_N}{2\eps} + O(1)
\end{align}
where $O(1)$ is bounded as $\eps\to 0$ independently of $j$.
Letting $j \to + \infty$ yields
\begin{align}
\label{formula from equation} \eps k_N \int_B u \; |x|^{-N+2\eps}
\le \omega_N \frac{H_N}{2\eps} + O(1),
\end{align}
showing that the integral on the left hand side is finite. On the
other hand, by \eqref{cotainf}
\begin{align}
\label{lower bound} \eps k_N \int_B u \; |x|^{-N+2\eps} \ge \eps
k_N \omega_N (1-\sigma) \int_0^1 \log(\frac{1}{r^4}) r^{-1+2\eps}
\, d r   = k_N \omega_N (1-\sigma) \frac{1}{\eps}.
\end{align}
Combining \eqref{formula from equation} and \eqref{lower bound} we
obtain
\begin{align*}
(1-\sigma) k_N \le \frac{H_N}2 + O(\eps).
\end{align*}
Letting $\eps \to 0 $ and then $\sigma \to 0$ we have
\begin{align*}
8 (N-2)(N-4) \le H_N = \frac{N^2 (N-4)^2}{16}.
\end{align*}
This is valid only if $N\ge 13$, a contradiction. \qed

\begin{rem}
The conclusion of Theorem~\ref{t1} can be obtained also from
Proposition~\ref{weakly singular}. However that proposition
depends crucially on the radial symmetry of the solutions, while
the argument in this section can be generalized to other domains.
\end{rem}

%%%%%%%%%%%%%%%%%%%%%%%%%%%%%%%%%%%%%%%%%%%%%%%%%%%%%%%%%%%%%%%%%%%%%%

\section{The extremal solution is singular in large
dimensions}\label{s-unbounded-1}

{In this section we take $a=b=0$ and prove Theorem~\ref{t3} for
$N\ge 32$.}

The idea for the proof of Theorem~\ref{t3} is to to estimate
accurately from above the function $\lambda^* e^{u^*}$, and to
deduce that the operator $\Delta^2 - \lambda^* e^{u^*}$ has a
strictly positive first eigenvalue (in the $H_0^2(B)$ sense).
Then, necessarily, $u^*$ is singular.

Upper bounds for both $\lambda^*$ and $u^*$ are obtained by
finding suitable sub and supersolutions. For example, if for some
$\lambda_1$ there exists a supersolution then $\lambda^* \ge
\lambda_1$. If for some $\lambda_2$ one can exhibit a stable
singular subsolution $u$, then $\lambda^* \le \lambda_2$.
Otherwise $\lambda_2 < \lambda^*$ and one can then prove that the
minimal solution $u_{\lambda_2}$ is above $u$, which is
impossible. The bound for $u^*$ also requires a stable singular
subsolution.

It turns out that in dimension $N\ge 32$ we can construct the
necessary subsolutions and verify their stability by hand. For
dimensions $13 \le N \le 31$ it seems difficult to find these
subsolutions explicitly. We adopt then an approach that involves a
computer assisted construction of subsolutions and verification of
the desired inequalities. We present this part in the next
section.

%We shall however begin by a result that allows us to consider only
%the case $b=0$ in Theorem~\ref{t3}.

\begin{lem} \label{bound-u*}
Assume $N\ge 13$. Then $ u^* \le \bar u =-4\log|x|$ in $B_1$.
\end{lem}
\proof{} Define $ \bar u(x) = -4 \log|x|$. Then $\bar u $
satisfies
\begin{equation*}
\left\{
\begin{aligned}
\Delta^2 \bar u &= 8 (N - 2) ( N - 4) e^{\bar u} && \hbox{in
$\R^N$}
\\
\bar u &= 0 && \hbox{on $\pa B_1$}
\\
\pd{\bar u}n &= -4 && \hbox{on $\pa B_1$}
\end{aligned}
\right.
\end{equation*}

Observe that since $\bar u$ is a supersolution to \eqref{main2}
{with $a=b=0$} we deduce immediately that $ \lambda^* \ge 8 (N -
2) ( N - 4)$.

In the case $ \la^* = 8 (N - 2) ( N - 4)$ we have $ u_\la \le \bar
u$ for all $0 \le \la <\la^*$ because $ \bar u$ is a
supersolution, and therefore $ u^* \le \bar u$ holds.
Alternatively, one can invoke Theorem 3 in \cite{aggm} to conclude
that we always have $ \la^* > 8 (N - 2) ( N - 4)$.

Suppose now that $ \la^* > 8 (N - 2) ( N - 4)$. We prove that $
u_\la \le \bar u$ for all $ 8 (N - 2) ( N - 4) < \la <\la^*$. Fix
such $\la$ and assume by contradiction that $ u_\la \le \bar u$ is
not true. Note that for $r<1$ and sufficiently close to 1 we have
$u_\lambda(r)<\bar u(r)$ because $u_\lambda'(1)=0$ while $\bar
u'(1)= -4$. Let
\begin{equation*}
R_1 = \inf \{ \, 0 \le R \le 1 \, \mid \, u_\la < \bar u \hbox{ in
} (R,1) \, \}.
\end{equation*}
Then $ 0 < R_1 < 1$, $u_\la(R_1) = \bar u(R_1)$ and $ u_\la'(R_1)
\le \bar u'(R_1)$. {So $ u_\la$ is a super-solution to the problem
\begin{equation}
\label{a2} \left\{
\begin{aligned}
\Delta^2 u &= 8 (N-2)(N-4)  e^u && \hbox{in $B_{R_1}$}
\\
u &= u_\la(R_1) && \hbox{on $\pa B_{R_1}$}
\\
\pd un &= u_\la'(R_1) && \hbox{on $\pa B_{R_1}$}
\end{aligned}
\right.
\end{equation}
while $\bar u$ is a subsolution to \eqref{a2}. Moreover it is
stable for this problem, since from Rellich's inequality
\eqref{hardy-rellich} and $ 8 (N-2)(N-4) \le N^2 (N-4)^2 / 16$ for
$N \ge 13$, we have
\begin{align*}
8 (N-2)(N-4) \int_{B_{R_1}} e^{\bar u} {\varphi^2} \le \frac{N^2
(N-4)^2}{16} \int_{\R^N} \frac{\varphi^2}{|x|^4} \le \int_{\R^N} (
\Delta \varphi)^2  \quad \forall \varphi \in C_0^\infty(B_{R_1}).
\end{align*}}
By Remark~\ref{comp-non-linear-2} we deduce $ \bar u \le u_\la$ in
$B_{R_1}$ which is impossible. \qed

%for the cases $13\le N\le 31$ and $N\ge 32$
%
%
%share the same basic idea. But we can ca
%
%
%separately

An upper bound for $\lambda^*$ is obtained by considering again a
stable, singular subsolution to the problem (with another
parameter, though):
\begin{lem}\label{L50}
For $N\ge 32$ we have
\begin{align} \label{100}
\la^* \le 8(N-2)(N-4)e^2.
\end{align}
\end{lem}
\proof{} Consider $w = 2 (1-r^2)$ and define
\begin{align*}
u= \bar u - w
\end{align*}
where $\bar u(x)= -4 \log |x|$. Then
\begin{align*}
\Delta^2 u = 8 (N-2)(N-4) \frac 1{r^4} = 8 (N-2)(N-4) e^{\bar u}
&= 8 (N-2)(N-4) e^{u + w}
\\
& \le 8 (N-2)(N-4)e^2 e^u.
\end{align*}
Also $u(1) = u'(1) =0$, so $u$ is a subsolution to \eqref{main2}
with parameter $\la_0 = 8 (N-2)(N-4) e^2$.

%There is an improved Hardy-Rellich inequality \cite{ggm} which
%states that if $ N \ge 5$ and $\Omega \subset \R^N$ is an open
%bounded set then
%\begin{equation
%\int_{\Omega} ( \Delta \varphi)^2 \ge \frac{N^2 (N-4)^2}{16}
%\int_{\Omega} \frac{\varphi^2}{|x|^4} + C_1 \int_{\Omega}
%\frac{\varphi^2}{|x|^2} + C_2 \int_{\Omega} \varphi^2 \quad
%\forall \varphi \in C_0^\infty(\Omega).
%\end{equation}
%where $C_1,C_2>0$. The constant $ \frac{N^2 (N-4)^2}{16}$ is the
%best possible.

For $N \ge 32$ we have $\la_0 \le N^2 (N-4)^2 /16$. Then by
\eqref{hardy-rellich} $u$ is a stable subsolution of \eqref{main2}
with $\lambda=\lambda_0$. If $\la^*
> \lambda_0 = 8(N-2)(N-4)e^2$ the minimal solution $u_{\lambda_0}$ to
\eqref{main2} with parameter $\lambda_0$ exists and is smooth.
From Lemma~\ref{comp-non-linear} we find $ u \le u_{\la_0}$ which
is impossible because $ u$ is singular and $u_{\la_0}$ is bounded.
Thus we have proved \eqref{100} for $N\ge 32$. \qed

\proof{ of Theorem~\ref{t3} in the case $N\ge 32$}

Combining Lemma~\ref{bound-u*} and \ref{L50} we have that if $ N
\ge 32$ then $ \lambda^* e^{u^*} \le r^{-4} \, 8 (N-2)(N-4)e^2 \le
r^{-4} N^2 (N-4)^2/ 16 $. This and \eqref{hardy-rellich} show that
\begin{align*}
\inf_{ \varphi \in C_0^\infty(B) } \frac{ \int_{B} (\Delta
\varphi)^2 - \la^* \int_{B} e^{u^*} \varphi^2}{\int_{B} \varphi^2}
> 0
\end{align*}
which is not possible if $u^*$ is bounded. \qed

%This was exactly the argument in the proof of Theorem~\ref{t3}
%where the stable subsolution was given by $\bar u - w$. We state
%the conditions needed in Lemma~\ref{L39} below. A lower bound for
%$\lambda^*$ is also need, and to obtain it one needs a
%supersolution. This is Lemma~\ref{L40}. Then finally we obtain an
%upper bound for $u^*$

\section{A computer assisted proof for dimensions $13\le N \le
31$}\label{s-unbounded-2}

{Throughout this section we assume $a=b=0$.} As was mentioned in
the previous section, the proof of Theorem~\ref{t3} relies on
precise estimates for $u^*$ and $\lambda^*$. We present first some
conditions under which it is possible to find these estimates.
Later we show how to meet such conditions with a computer assisted
verification.

The first lemma is analogous to Lemma~\ref{L50}.

\begin{lem}\label{L39}
Suppose there exist $\eps>0$, $\lambda>0$ and a radial function $u
\in H^2(B) \cap W^{4,\infty}_{loc}(B\setminus\{0\})$ such that
\begin{align}
& \nonumber \Delta^2 u \le \lambda e^u \quad \hbox{ for all }
0<r<1
\\
& \nonumber|u(1)| \le \eps , \quad \left| \pd{u}{n}(1) \right| \le
\eps
\\
& \nonumber u \not\in L^\infty(B)
\\
& \label{700} \lambda e^{\eps} \int_B e^u \varphi^2 \le \int_B
(\Delta \varphi)^2 \quad \hbox{ for all } \varphi \in
C_0^\infty(B).
\end{align}
Then
\begin{align*}%\label{lambda-max}
\lambda^* \le \lambda e^{2 \eps}.
\end{align*}
\end{lem}
\proof{} Let
\begin{align}\label{710}
\psi(r) = \eps r^2 - 2 \eps
\end{align}
so that
\begin{align*}
\Delta^2 \psi \equiv 0 , \quad \psi(1) = - \eps , \quad \psi'(1) =
2 \eps
\end{align*}
and
\begin{align*}
-2\eps \le\psi(r) \le -\eps \quad \hbox{ for all } 0\le r \le 1.
\end{align*}
It follows that
\begin{align*}
\Delta^2 (u + \psi) &\le \lambda e^u = \lambda e^{-\psi}
e^{u+\psi} \le \lambda e^{2 \eps} e^{u+\psi} .
\end{align*}
On the boundary we have $u(1) + \psi(1)\le 0 $, $u'(1) + \psi'(1)
\ge 0$. Thus $u+\psi$ is a singular subsolution to the equation
with parameter $ \lambda e^{2 \eps} $. Moreover, since $\psi \le
-\eps $ we have $\lambda e^{2\eps} e^{u + \psi} \le \lambda e^\eps
e^u $ and hence, from \eqref{700} we see that $u+ \psi$ is stable
for the problem with parameter $\lambda e^{2\eps} $. If $\lambda
e^{2\eps} < \lambda^*$ then the minimal solution associated to the
parameter $ \lambda e^{2\eps}$ would be above $u+ \psi$, which is
impossible because $u$ is singular. \qed

\begin{lem}\label{L40}
Suppose we can find $\eps>0$, $\lambda>0$ and $u \in H^2(B) \cap
W^{4,\infty}_{loc}(B\setminus\{0\})$ such that
\begin{align}
& \nonumber \Delta^2 u \ge \lambda e^u \quad \hbox{ for all }
0<r<1
\\
& \nonumber|u(1)| \le \eps , \quad \left| \pd{u}{n}(1) \right| \le
\eps .
\end{align}
Then
\begin{align*}%\label{lambda-min}
\lambda e^{- 2 \eps} \le \lambda^*.
\end{align*}
\end{lem}
\proof{} Let $\psi$ be given by \eqref{710}. Then $u - \psi$ is a
supersolution to the problem with parameter $\lambda e^{-2 \eps}$.
\qed

The next result is the main tool to guarantee that $u^*$ is
singular. The proof, as in Lemma~\ref{bound-u*}, is based on an
upper estimate of $u^*$ by a stable singular subsolution.
\begin{lem}\label{L51}
Suppose there exist $\eps_0,\eps>0$, $\lambda_a > 0$ and a radial
function $u \in H^2(B) \cap W^{4,\infty}_{loc}(B\setminus\{0\})$
such that
\begin{align}
& \label{1000} \Delta^2 u \le ( \lambda_a+\eps_0) e^u \quad \hbox{
for all } 0<r<1
\\
& \label{1010} \Delta^2 u \ge ( \lambda_a - \eps_0) e^u \quad
\hbox{ for all } 0<r<1
\\
& |u(1)| \le \eps , \quad \left| \pd{u}{n}(1) \right| \le \eps
\\
& \label{1040} u \not\in L^\infty(B)
\\
& \label{460} \beta_0 \int_B e^u \varphi^2 \le \int_B (\Delta
\varphi)^2 \quad \hbox{ for all } \varphi \in C_0^\infty(B),
\end{align}
where
\begin{align}\label{def-beta0}
\beta_0 = \frac{(\lambda_a + \eps_0)^3 }{(\lambda_a - \eps_0)^2 }
e^{9 \eps} .
\end{align}
Then $u^*$ is singular and
\begin{align}\label{1320}
(\lambda_a - \eps_0) e^{-2\eps} \le \lambda^* \le (\lambda_a +
\eps_0) e^{2\eps}.
\end{align}
\end{lem}
\proof{} By Lemmas~\ref{L39} and \ref{L40} we have \eqref{1320}.
Let
\begin{align*}
\delta = \log\left( \frac{\lambda_a + \eps_0}{\lambda_a -\eps_0}
\right) + 3 \eps.
\end{align*}
and define
\begin{align*}
\varphi (r) = -\frac{\delta}{4} r^4 + 2 \delta.
\end{align*}
We claim that
\begin{align}\label{440}
u^* \le u + \varphi \quad \hbox{ in } B_1.
\end{align}
To prove this, we shall show that for
%$(\lambda_a-\eps) e^{-\delta} < \lambda < \lambda^*$
$\lambda < \lambda^*$
\begin{align}\label{430}
u_\lambda \le u + \varphi \quad \hbox{ in } B_1 .
\end{align}
Indeed, we have
\begin{align*}
& \Delta^2 \varphi = - \delta 2 N(N+2)
\\
& \varphi(r) \ge \delta \quad \hbox{ for all } 0 \le r \le 1
\\
& \varphi(1) \ge \delta \ge \eps, \quad \varphi'(1) = - \delta \le
- \eps
\end{align*}
and therefore
\begin{align}
\nonumber \Delta^2(u+\varphi) & \le (\lambda_a+\eps_0) e^u +
\Delta^2 \varphi \le (\lambda_a+\eps_0) e^u = (\lambda_a+\eps_0)
e^{-\varphi} e^{u + \varphi}
\\
\label{420} & \le (\lambda_a+\eps_0) e^{-\delta} e^{u + \varphi}.
\end{align}
%From Lemma~\ref{L39} we deduce $\lambda^* \le (\lambda_a + \eps) e^{2\eps}$ and
By \eqref{1320} and the choice of $\delta$
\begin{align}\label{450}
(\lambda_a + \eps_0) e^{-\delta} = (\lambda_a - \eps_0 ) e^{-3
\eps} < \lambda^* .
\end{align}
To prove \eqref{430} it suffices to consider $\lambda$ in the
interval $(\lambda_a-\eps_0 ) e^{-3 \eps} < \lambda < \lambda^*$.
Fix such $\lambda$ and assume that \eqref{430} is not true. Write
\begin{align*}
\bar u = u + \varphi
\end{align*}
and let
\begin{equation*}
{R_1 = \sup \{ \, 0 \le R \le 1 \, |\, u_\la(R) = \bar u(R) \, \}.
}
\end{equation*}
%{\color{red} I think this is the same as the formula that was here
%before
%\begin{equation*}
%R_1 = \max \{ \, 0 \le R \le 1 \, |\, u_\la(R) = \bar u(R) \, \}.
%\end{equation*}
%No? }

Then $ 0 < R_1 < 1$ and $ u_\la(R_1) = \bar u(R_1)$. Since $
u_\la'(1)=0$ and $\bar u '(1) < 0$ we must have $ u_\la'(R_1) \le
\bar u'(R_1)$. Then $ u_\la$ is a solution to the problem
\begin{equation*}
\left\{
\begin{aligned}
\Delta^2 u &= \la e^u && \hbox{in $B_{R_1}$}
\\
u &= u_\la(R_1) && \hbox{on $\pa B_{R_1}$}
\\
\pd un &= u_\la'(R_1) && \hbox{on $\pa B_{R_1}$}
\end{aligned}
\right.
\end{equation*}
while, thanks to \eqref{420} and \eqref{450}, $\bar u$ is a
subsolution to the same problem. Moreover $\bar u$ is stable
thanks to \eqref{460} since, by Lemma~\ref{L39},
\begin{align}\label{470}
\lambda < \lambda^* \le (\lambda_a + \eps_0 ) e^{2\eps}
\end{align}
and hence
\begin{align*}
\lambda e^{\bar u} \le (\lambda_a + \eps_0) e^{2\eps} e^{2\delta }
e^u \le \beta_0 e^u.
\end{align*}
We deduce $ \bar u \le u_\la$ in $B_{R_1}$ which is impossible,
since $\bar u$ is singular while $u_\lambda$ is smooth. This
establishes \eqref{440}.

From \eqref{440} and \eqref{470} we have
\begin{align*}
\lambda^* e^{u^*} \le \beta_0 e^{-\eps} e^u
\end{align*}
and therefore
\begin{align*}
\inf_{\varphi \in C_0^\infty(B)} \frac{\int_B ( \Delta\varphi)^2 -
\lambda^* e^{u^*} \varphi^2 }{ \int_B \varphi^2 } >0.
\end{align*}
This is not possible if $u^* $ is a smooth solution. \qed

For each dimension $13\le N \le 31$ we construct $u$ satisfying
\eqref{1000} to \eqref{460} of the form
\begin{align}\label{def-u-num}
u(r) = \begin{cases} - 4 \log r + {\log\left(
\frac{8(N-2)(N-4)}{\lambda} \right)} & \hbox{ for } 0<r<r_0
\\
\tilde u(r) & \hbox{ for } r_0 \le r \le 1
\end{cases}
\end{align}
where $\tilde u$ is explicitly given. Thus $u$ satisfies
\eqref{1040} automatically.

Numerically it is better to work with the change of variables
\begin{align*}
w(s) = u(e^s) + 4 s, \quad -\infty<s<0
\end{align*}
which transforms the equation $\Delta^2 u = \lambda e^u$ into
\begin{align*}
L w + 8(N-2)(N-4) = \lambda e^w , \quad -\infty<s<0
\end{align*}
where
\begin{align*}
L w = \frac{d^4 w}{d s^4} + 2(N-4) \frac{d^3 w}{d s^3} + (N^2-10 N
+ 20 ) \frac{d^2 w}{d s^2} - 2 (N-2)(N-4) \frac{d w}{d s}
\end{align*}
The boundary conditions $u(1)=0$, $u'(1)=0$ then yield
\begin{align*}
w(0)=0, \quad w'(0)=4.
\end{align*}
Regarding the behavior of $w$ as $s\to -\infty$ observe that
$$
u(r) = - 4 \log r + { \log\left( \frac{8(N-2)(N-4)}{ \lambda}
\right)} \quad \hbox{for $r< r_0$ }
$$
if and only if
\begin{align*}
w(s) = \log \frac{8(N-2)(N-4)}{\lambda} \quad \hbox{for all } s <
\log r_0.
\end{align*}

The steps we perform are the following:
\medskip

1) We fix $x_0<0$ and using numerical software we follow a branch
of solutions to
\begin{align*}
\left\{
\begin{aligned}
L \hat w + 8(N-2)(N-4) &= \lambda e^{\hat w} , \quad x_0<s<0
\\
\hat w(0) &=0 , \quad \hat w'(0) = t
\\
\hat w(x_0) &= \log \frac{8(N-2)(N-4)}{\lambda}, \quad \frac{d^2
\hat w}{d s^2}(x_0)=0 ,  \quad \frac{d^3 \hat w}{d s^3}(x_0)=0
\end{aligned}
\right.
\end{align*}
as $t$ increases from 0 to 4. The numerical solution $ ( \hat w$,
$\hat \lambda)$ we are interested in corresponds to the case
$t=4$. The five boundary conditions are due to the fact that we
are solving a fourth order equation with an unknown parameter
$\lambda$.
\medskip

2) Based on $\hat w$, $\hat \lambda$ we construct a $C^3$ function
$w$ which is constant for $s\le x_0$ and piecewise polynomial for
$x_0\le s \le 0$. More precisely, we first divide the interval
$[x_0,0]$ in smaller intervals of length $h$. Then we generate a
cubic spline approximation $g_{fl} $ with floating point
coefficients of $\frac{d^4 \hat w}{ds^4}$. From $g_{fl}$ we
generate a piecewise cubic polynomial $g_{ra}$ which uses rational
coefficients and we integrate it 4 times to obtain $w$, where the
constants of integration are such that $\frac{d^j w}{d
s^j}(x_0)=0$, $1\le j \le 3$ and $w(x_0)$ is a rational
approximation of $\log( 8(N-2)(N-4) / \lambda)$. Thus $w$ is a
piecewise polynomial function that in each interval is of degree 7
with rational coefficients, and which is globally $C^3$. We also
let $\lambda$ be a rational approximation of $\hat \lambda$. With
these choices note that $L w + 8(N-2)(N-4) - \lambda e^w$ is a
small constant (not necessarily zero) for $s\le x_0$.
\medskip

3) The condition \eqref{1000} and \eqref{1010} we need to check
for $u$ are equivalent to the following inequalities for $w$
\begin{align}
\label{1200} L w + 8(N-2)(N-4) - (\lambda+\epsilon_0) e^w \le 0,
\quad -\infty<s<0
\\
\label{1210} L w + 8(N-2)(N-4) - (\lambda - \epsilon_0) e^w \ge 0,
\quad -\infty<s<0 .
\end{align}
Using a program in Maple we verify that $w$ satisfies \eqref{1200}
and \eqref{1210}. This is done evaluating a second order Taylor
approximation of $L w + 8(N-2)(N-4) - (\lambda+\epsilon) e^w $ at
sufficiently close mesh points. All arithmetic computations are
done with rational numbers, thus obtaining exact results. The
exponential function is approximated by a Taylor polynomial of
degree 14 and the difference with the real value is controlled. {

More precisely, we write
\begin{align*}
f(s) &= L w + 8(N-2)(N-4) - (\lambda+\epsilon_0) e^w,
\\
\tilde f(s) &= L w + 8(N-2)(N-4) - (\lambda+\epsilon_0) T(w),
\end{align*}
where $T$ is the Taylor polynomial of order 14 of the exponential
function around 0. Applying Taylor's formula to $f$ at $y_j$, we
have for $s\in[y_j,y_{j+h}]$,
\begin{align*}
f(s) &\le f(y_j) + |f'(y_j)|h + \frac12 Mh^2
\\
&\le \tilde f(y_j) + |\tilde f'(y_j)|h + \frac12 Mh^2 +
|f(y_j)-\tilde f(y_j)| + |f'(y_j)-\tilde f'(y_j)|h
\\
&\le \tilde f(y_j) + |\tilde f'(y_j)| {h} + \frac{1}{2} M {h}^2 +
E_1 + E_2h,
\end{align*}
where
\begin{align*}
&\hbox{$M$ is a bound for $|f''|$ in $[y_j,y_j+h]$}
\\
&\hbox{$E_1$ is such that $(\lambda+\epsilon_0 )|e^w - T(w)|\le
E_1$ in $[y_j,y_j+h]$}
\\
&\hbox{$E_2$ is such that $(\lambda+\epsilon_0 )|( e^w - T'(w) )
w' |\le E_2$ in $[y_j,y_j+h]$.}
\end{align*}
So, inequality \eqref{1200} will be verified on each interval
$[y_j, y_j+h]$ where $w$ is a polynomial as soon as
\begin{align}\label{numerical-ineq}
\tilde f(y_j) + | \tilde f'(y_j)| {h} + \frac{1}{2} M {h}^2 + E_1
+ E_2h \le 0.
\end{align}
When more accuracy is desired, instead of \eqref{numerical-ineq}
one can verify that
$$
\tilde f(x_i) + | \tilde f'(x_i)| \frac{h}{m} + \frac{1}{2} M
(\frac{h}{m})^2 + E_1 + E_2\frac h m \le 0,
$$
where $(x_i)_{i=1\dots m+1}$ are $m+1$ equally spaced points in
$[y_j,y_j+h]$.

We obtain exact values for the upper bounds $M, E_1, E_2$ as
follows. First note that $f'' = L w''-(\lambda+\eps_0) e^w (
(w')^2 + w'')$. On $[y_j,y_j+h]$, we have $w(s) = \sum_{i=0}^7 a_i
(s-y_j)^i$ and we estimate $|w(s)| \le \sum_{i=0}^7 |a_i| h^i $
for $s\in [y_j,y_j+h]$. Similarly,
\begin{align}\label{bound derivative}
\left| \frac{d^\ell w}{d s^\ell}(s) \right| \le \sum_{i=\ell}^7 i
(i-1) \ldots (i-\ell+1) |a_i| h^{i-\ell} \quad \hbox{for all }
s\in [y_j,y_j+h]
\end{align}
The exponential is estimated by $e^w \le e^1 \le 3$, since our
numerical data satisfies the rough bounds $-3/2\le w \le 1$. Using
this information and \eqref{bound derivative} yields a rational
upper bound $M$. $E_1$ is estimated using Taylor's formula~:
\begin{align*}
E_1 = (\lambda+\eps_0) \frac{(3/2)^{15} }{15!}.
\end{align*}
Similarily, $E_2 = (\lambda+\eps_0) \frac{(3/2)^{14} }{14!} B_1$
where $B_1$ is the right hand side of \eqref{bound derivative}
when $\ell=1$. }
%{\color{red} I just realized that the Maple program is not really as
%described above. In the current version we are using $e^w \le 1$.
%This should be corrected, and the we should run the tests again.
%Maybe the bound $e^w \le e^1 \le 3$ is bad ... }
\medskip

4) We show that the operator $\Delta^2 - \beta e^{u}$ where
$u(r)=w(\log r)-4\log r$, satisfies condition \eqref{460} for some
$\beta \ge \beta_0$ where $\beta_0$ is given by \eqref{def-beta0}.
In dimension $N\ge 13$ the operator $\Delta^2 - \beta e^{u}$ has
indeed a positive eigenfunction in $H_0^2(B)$ with finite
eigenvalue if $\beta$ is not too large. The reason is that near
the origin
\begin{align*}
\beta e^u = \frac{c}{|x|^4},
\end{align*}
where $c$ is a number close to $ 8 (N-2)(N-4) \beta / \lambda $.
If $\beta$ is not too large compared to $\lambda$ then $c < N^2
(N-4)^2 /16$ and hence, using \eqref{hardy-rellich}, $\Delta^2 -
\beta e^u$ is coercive in $H_0^2(B_{r_0})$ (this holds under even
weaker conditions, see \cite{davila-dupaigne}). It follows that
there exists a first eigenfunction $\varphi_1 \in H_0^2(B) $ for
the operator $\Delta^2 - \beta e^u$ with a finite first eigenvalue
$\mu_1$, that is
\begin{align*}
&\Delta^2 \varphi_1 - \beta e^u \varphi_1 = \mu_1 \varphi_1 \quad
\hbox{in $B$}
\\
& \varphi_1 >0 \quad \hbox{in $B$}
\\
&\varphi_1 \in H_0^2(B).
\end{align*}
Moreover $\mu_1$ can be characterized as
\begin{align*}
\mu_1 = \inf_{\varphi \in C_0^\infty(B)} \frac{\int_B (\Delta
\varphi)^2 - \beta e^u \varphi^2}{\int_B \varphi^2}
\end{align*}
and is the smallest number for which a positive eigenfunction in
$H_0^2(\Omega)$ exists.

Thus to prove that \eqref{460} holds it suffices to verify that
$\mu_1\ge 0$ and for this it is enough to show the existence of a
nonnegative $\varphi \in H_0^2(B)$, $\varphi \not \equiv 0$ such
that
\begin{align}\label{1300}
\left\{
\begin{aligned}
\Delta^2 \varphi - \beta e^u \varphi &\ge 0 && \hbox{in $B$}
\\
\varphi &= 0 && \hbox{on $\partial B$}
\\
\pd{\varphi}{n} &\le 0 && \hbox{on $\partial B$.}
\end{aligned}
\right.
\end{align}
Indeed, multiplication of \eqref{1300} by $\varphi_1$ and
integration by parts yields
\begin{align*}
\mu_1 \int_B \varphi \varphi_1 + \int_{\partial B} \pd{\varphi}{n}
\Delta \varphi_1 \ge 0.
\end{align*}
But $\Delta \varphi_1 \ge 0$ on $\partial B$ and thus $\mu_1 \ge
0$. To achieve \eqref{1300} we again change variables and define
\begin{align*}
\phi(s) = \varphi(e^s) \quad -\infty<s\le 0.
\end{align*}
Then we have to find $\phi \ge 0$, $\phi\not\equiv 0$ satisfying
\begin{align}\label{1310}
\left\{
\begin{aligned}
 L \phi - \beta e^ w \phi &\ge 0 \quad \hbox{in $-\infty<s\le 0$}
\\
\phi(0) &=0
\\
\phi'(0) & \le 0.
\end{aligned}
\right.
\end{align}
Regarding the behavior as $s \to -\infty$, we note that $w$ is
constant for $-\infty<s<x_0$, and therefore, if
\begin{align*}
L\phi - \beta e^w \phi \equiv 0 \quad -\infty<s\le x_0
\end{align*}
then $\phi$ is a linear combination of exponential functions
$e^{-\alpha s}$ where $\alpha$ must be a solution to
\begin{align*}
\alpha^4 - 2(N-4) \alpha^3 + (N^2-10 N + 20 ) \alpha^2 + 2
(N-2)(N-4) \alpha = \beta e^{w(x_0)}
\end{align*}
where $ \beta e^{w(x_0)}$ is close to $8(N-2) (N-4)
\beta/\lambda$. If $N\ge 13$ the polynomial
\begin{align*}
\alpha^4 - 2(N-4) \alpha^3 + (N^2-10 N + 20 ) \alpha^2 + 2
(N-2)(N-4) \alpha - 8(N-2)(N-4)
\end{align*}
has 4 distinct real roots, while if $N\le 12$ there are 2 real
roots and 2 complex conjugate. If $N\ge 13$ there is exactly one
root in the interval $(0, (N-4)/2 )$, 2 roots greater than
$(N-4)/2$ and one negative. We know that $\varphi(r) = \phi(\log
r) \sim r^{-\alpha}$ is in $H^2$, which forces $\alpha < (N-4)/2$.
It follows that for $s <x_0$, $\phi $ is a combination of
$e^{-\alpha_0 s}$, $e^{-\alpha_1 s}$ where $\alpha_0>0$,
$\alpha_1<0$ are the 2 roots smaller than $\alpha < (N-4)/2$. For
simplicity, however we will look for $\phi$ such that $\phi(s) = C
e^{-\alpha_0s} $ for $s<x_0$ where $C>0$ is a constant. This
restriction will mean that we will not be able to impose
$\phi'(0)=0$ at the end. This is not a problem because
$\phi'(0)\le 0$.

Notice that we only need the inequality in \eqref{1310} and hence
we need to choose $\alpha \in (0,N-4/2)$ such that
\begin{align*}
\alpha^4 - 2(N-4) \alpha^3 + (N^2-10 N + 20 ) \alpha^2 + 2
(N-2)(N-4) \alpha \ge \beta e^{w(x_0)}.
\end{align*}
The precise choice we employed in each dimension is in a summary
table at the end of this section.

To find a suitable function $\phi$ with the behavior $\phi(s) = C
e^{-\alpha} $ for $s<x_0$ we set $\phi = \psi e^{-\alpha s} $ and
solve the following equation
\begin{align*}
T_\alpha \psi - \beta e^w \psi = f
\end{align*}
where the operator $T_\alpha$ is given by
\begin{align*}
T_\alpha \psi &= \frac{d^4 \psi}{d s^4} + ( - 4 \alpha + 2(N-4) )
\frac{d^3 \psi}{d s^3} + ( 6 \alpha^2 - 6 \alpha (N-4) + N^2-10 N
+ 20 ) \frac{d^2 w}{d s^2}
\\
& \quad + ( - 4 \alpha^3 + 6 \alpha^2 (N-4) -2\alpha (N^2-10 N +
20 ) - 2 (N-2)(N-4) ) \frac{d \psi}{d s}
\\
&\quad + ( \alpha^4 - 2 \alpha^3 (N-4) + \alpha^2 (N^2 -10 N + 20)
+ 2 \alpha (N-2) (N-4) ) \psi
\end{align*}
and $f$ is some smooth function such that $f\ge 0$,
$f\not\equiv0$. Actually we choose $\bar \beta > \beta_0$ (where
$\beta_0$ is given in \eqref{def-beta0}) find $\bar \alpha$
satisfying approximately
\begin{align*}
\bar \alpha^4 - 2(N-4) \bar \alpha^3 + (N^2-10 N + 20 ) \bar
\alpha^2 + 2 (N-2)(N-4) \bar \alpha = \bar \beta e^{w(x_0)}.
\end{align*}
We solve numerically
\begin{align*}
T_{\bar \alpha} \hat\psi - \bar \beta e^w \hat\psi &= f \quad
x_0<s<0
\\
\hat\psi(x_0) & =1, \quad \hat\psi''(x_0) = 0, \quad
\hat\psi'''(x_0) = 0
\\
\hat\psi(0) &= 0
\end{align*}
Using the same strategy as in 2) from the numerical approximation
of $\frac{d^4 \hat\psi}{d s^4}$ we compute a piecewise polynomial
$\psi$ of degree 7, which is globally $C^3$ and constant for $s\le
x_0$. The constant $\psi(x_0)$ is chosen so that $\psi(0)=0$. We
use then Maple to verify the following inequalities
\begin{align*}
\psi & \ge 0 \quad x_0 \le s \le 0
\\
T_\alpha \psi - \beta e^w \psi & \ge 0 \quad x_0 \le s \le 0
\\
\psi'(0) & \le 0
\end{align*}
where $\beta_0 < \beta< \bar \beta$ and $0 < \alpha < (N-4)/2$ are
suitably chosen.

\bigskip
{\cb At the URLs: \smallskip

\noindent
{\tt http://www.lamfa.u-picardie.fr/dupaigne/ \\
http://www.ime.unicamp.br/\~\ \!\!\!\!\!\! msm/} \smallskip

\noindent we have provided the data of the functions $w$ and $\psi$
defined as piecewise polynomials of degree 7 in $[x_0, 0 ]$ with
rational coefficients for each dimension in $13\leq  N\leq 31.$ We
also give a rational approximation of the constants involved in the
corresponding problems.}

{\cb We use Maple to verify that $w$ and $\psi$ (with suitable
extensions) are $C^3$ global functions and satisfy the corresponding
inequalities, using only its capability to operate on arbitrary
rational numbers. These operations are exact and are limited only by
the memory of the computer and clearly slower than floating point
operations. We chose Maple since it is a widely used software, but
the reader can check the validity of our results with any other
software (see e.g. the open-source solution pari/gp).

The tests were conducted using Maple 9. }

\bigskip

Summary of parameters and results

\medskip

\begin{center}

\begin{supertabular}{ccccccc}
%\begin{tabular}{ccccccc}
$N$ & $\lambda$ & $\epsilon_0$ & $\epsilon$ & $\bar\beta$ &
$\beta$ & $\alpha$
\\
\hline 13 & 2438.6 & 1 & $5\cdot 10^{-7}$ & 2550 & 2500 & 3.9
\\
14 & 2911.2 & 1 & $3 \cdot 10^{-6}$ & 3100 & 3000 & 3.4
\\
15 & 3423.8 & 1 & $3\cdot 10^{-6}$ & 3600 & 3500 & 3.1
\\
16 & 3976.4 & 1 & $1 \cdot 10^{-5}$ & 4100 & 4000 & 3.0
\\
17 & 4568.8 & 1 & $2\cdot 10^{-4}$ & 4800 & 4600 & 3.0
\\
18 & 5201.1 & 2 & $2\cdot 10^{-4}$ & 5400 & 5300 & 2.7
\\
19 & 5873.2 & 2 & $2\cdot 10^{-4}$ & 6100 & 6000 & 2.7
\\
20 & 6585.1 & 3 & $7\cdot 10^{-4}$ & 7000 & 6800 & 2.7
\\
21 & 7336.7 & 3 & $7\cdot 10^{-4}$ & 7700 & 7500 & 2.6
\\
22 & 8128.1 & 4 & $1\cdot 10^{-3}$ & 8600 & 8400 & 2.6
\\
23 & 8959.1 & 4 & $1\cdot 10^{-3}$ & 9400 & 9200 & 2.5
\\
24 & 9829.8 & 4 & $1\cdot 10^{-3}$ & 10400 & 10200 & 2.5
\\
25 & 10740.1 & 4 & $1\cdot 10^{-3}$ & 11400 & 11200 & 2.5
\\
26 & 11690.1 & 6 & $2\cdot 10^{-3}$ & 12400 & 12200 & 2.5
\\
27 & 12679.7 & 7 & $2\cdot 10^{-3}$ & 13400 & 13200 & 2.4
\\
28 & 13709.0 & 7 & $2\cdot 10^{-3}$ & 14500 & 14300 & 2.4
\\
29 & 14777.8 & 7 & $2\cdot 10^{-3}$ & 15400 & 15200 & 2.4
\\
30 & 15886.2 & 8 & $2\cdot 10^{-3}$ & 16600 & 16400 & 2.4
\\
31 & 17034.3 & 10 & $2\cdot 10^{-3}$ & 17600 & 17500 & 2.3
\\
\end{supertabular}
%\end{tabular}

\end{center}

\begin{rem}
1) Although we work with $\lambda$ rational, in the table above we
prefer to display a decimal approximation of $\lambda$.

2) In the previous table we selected a ``large'' value of
$\epsilon_0$ in order to have a fast verification with Maple. By
requiring more accuracy in the numerical calculations, using a
smaller value of $\epsilon_0$ and using more subintervals to
verify the inequalities in the Maple program it is possible to
obtain better estimates of $\lambda^*$. For instance, using
formulas \eqref{1320}, we obtained

\medskip

\begin{tabular}{ccccccccc}
$N$ & $\lambda$ & $\epsilon_0$ & $\epsilon$ & $\lambda^*_{min}$ &
$\lambda^*_{max}$ & $\bar\beta$ & $\beta$ & $\alpha$
\\
\hline
 13 & 2438.589 & 0.003 & $5\cdot 10^{-7}$ & 2438.583 & 2438.595 & 2550 & 2510 & 3.9
\\
14 & 2911.194 & 0.003 & $5\cdot 10^{-7}$ & 2911.188 & 2911.200 &
3100 & 3000 & 3.4
\end{tabular}

\medskip
The verification above, however, required to check 1500
subintervals of each of the 4500 intervals of length 0.002, which
amounts to substantial computer time.

\end{rem}

\section{Proof of Proposition \ref{p1.7}} \label{s-last}
Throughout this section, we restrict, as permitted, to the case
$a=0$.

a)Let $u$ denote the extremal solution of \eqref{main2} with
homogeneous Dirichlet boundary condition $a=b=0$. We extend $u$ on
its maximal interval of existence $(0,\bar R)$.
\begin{lem}\label{lemma71}
$\bar R<\infty$ and $u(r)\sim \log(\bar R-r)^{-4}$ for $r\sim\bar
R$.
\end{lem}

\proof{} The fact that $\bar R<\infty$ can be readily deduced from
Section 2 of \cite{agg}. We present an alternative (and more
quantitative) argument. We first observe that
\begin{equation}\label{diff1}
u'' - \frac1r u' >0\qquad\text{for all $r\in[1,\bar R)$.}
\end{equation}
Integrate indeed \eqref{main2} over a ball of radius $r$ to
conclude that
\begin{equation}\label{intbr}
0<\lambda\int_{B_r} e^u = \int_{\partial B_r} \pd{}{r}\Delta u =
\omega_N r^{N-1}\left(u''' +\frac{N-1}r (u''-\frac1r u')\right)
\end{equation}
If $r=1$, since $u$ is nonnegative in $(0,1)$ and $u(1)=u'(1)=0$,
we must have $u''(1)\ge0$. In fact, $u''(1)>0$. Otherwise, we
would have $u''(1)=0$ and $u'''(1)>0$ by \eqref{intbr},
contradicting $u>0$ in $(0,1)$. So, we may define
$$
R=\sup\{r>1 :u''(t) - \frac1t u'(t) >0\;\;\text{for all
$t\in[1,r)$}\}
$$
and we just need to prove that $R=\bar R$. Assume this is not the
case, then $u''(R)-\frac1R u'(R)=0$ and $u'''(R)=\left(u''-\frac1R
u'\right)'(R)\le0$. This contradicts \eqref{intbr} and we have
just proved \eqref{diff1}. In particular, we see that $u$ is
convex increasing on $(1,\bar R)$.

\

Since $u$ is radial, \eqref{main2} reduces to
\begin{equation}
\label{ode} u^{(4)} + \frac{2(N-1)}{r}u'''
+\frac{(N-1)(N-3)}{r^2}u'' -\frac{(N-1)(N-3)}{r^3}u' = \lambda
e^u.
\end{equation}
Multiply by $u'$ :
$$
u^{(4)}u' + \frac{2(N-1)}{r}u'''u' +\frac{(N-1)(N-3)}{r^2}u''u'
-\frac{(N-1)(N-3)}{r^3} {(u')^2} = \lambda (e^u)',
$$
which we rewrite as
\begin{align*}
\left[(u'''u')'-u'''u''\right] + 2(N-1)\left[\left(\frac1r
u''u'\right)' - u''\left(\frac1r u'\right)'\right] &\\
+(N-1)(N-3)\left(\frac{(u')^2}{2r^2}\right)' &=\lambda (e^u)'.
\end{align*}
By \eqref{diff1}, it follows that for $r\in [1,\bar R)$,
$$
\left[(u'''u')'-u'''u''\right] + 2(N-1)\left(\frac1r u''u'\right)'
+ (N-1)(N-3)\left(\frac{(u')^2}{2r^2}\right)'\ge\lambda (e^u)'.
$$
Integrating, we obtain for some constant $A$
$$
u'''u' - \frac{(u'')^2}{2} + 2(N-1)\frac1r u''u' +
\frac{(N-1)(N-3)}{2}\frac{(u')^2}{r^2}\ge\lambda e^u - A.
$$
We multiply again by $u'$ :
\begin{equation}\label{diff2}
\begin{aligned}
\left[(u''(u')^2)' - u''\left((u')^2\right)'\right] -\frac12
(u'')^2u'+ 2(N-1)\frac1r u''(u')^2&\\
+\frac{(N-1)(N-3)}{2}\frac1{r^2}(u')^3 &\ge(\lambda e^u -Au)'.
\end{aligned}
\end{equation}
We deduce from \eqref{diff1} that
$$\frac1r u''(u')^2 =
\frac12\left(\frac1r (u')^3\right)' - \frac12(u')^2(\frac1r
u')'\le\frac12 \left(\frac1r (u')^3\right)' \quad\text{and}$$
$$
\frac1{r^2}(u')^3\le \frac1r (u')^2u'' \le \frac12 \left(\frac1r
(u')^3\right)'.
$$
Using this information in \eqref{diff2}, dropping nonpositive
terms and integrating, we obtain for some constant $B$,
$$
u''(u')^2 + \frac{(N^2-1)}4\frac1r (u')^3\ge \lambda e^u -Au -B
$$
Applying \eqref{diff1} again, it follows that for
$C=\frac{N^2-1}{4}+1$
$$
C u''(u')^2 \ge \lambda e^u -Au -B
$$
which after multiplication by $u'$ and integration provides
positive constants $c,C$ such that
$$
(u')^4\ge c(e^u -Au^2 -Bu -C).
$$
At this point, we observe that since $u$ is convex and increasing,
$u$ converges to $+\infty$ as $r$ approaches $\bar R$. Hence, for
$r$ close enough to $\bar R$ and for $c>0$ perhaps smaller,
$$
u'\ge c\;e^{u/4}.
$$
By Gronwall's lemma, $\bar R$ is finite and
$$
u\le-4\log(\bar R -r)+ C \mbox{for} \:r \mbox{close to}\:\bar R.
$$

\

It remains to prove that $u\ge-4\log(\bar R -r)- C$. This time, we
rewrite \eqref{main2} as
$$
\left[r^{N-1}(\Delta u)'\right]' = \lambda r^{N-1}e^u.
$$
We multiply by $r^{N-1}(\Delta u)'$ and obtain :
$$
\frac12 \left[r^{2N-2}((\Delta u)')^2\right]' = \lambda
r^{2N-2}e^u(\Delta u)' \le C e^u(\Delta u)'\le C(e^u\Delta u)'.
$$
Hence, for $r$ close to $\bar R$ and $C$ perhaps larger,
$$
((\Delta u)')^2\le Ce^u\Delta u
$$
and so
$$
\sqrt{\Delta u}(\Delta u)'\le Ce^{u/2}\Delta u\le C'e^{u/2}u''\le
C'(e^{u/2}u')',
$$
where we have used \eqref{diff1}. Integrate to conclude that
$$
(\Delta u)^{3/2}\le C e^{u/2}u'.
$$
Solving for $\Delta u$ and multiplying by $(u')^{1/3}$, we obtain
in particular that
$$
(u')^{1/3}u''\le Ce^{u/3}u'.
$$
Integrating again, it follows that $(u')^{4/3}\le Ce^{u/3}$, i.e.
$$
u'\le Ce^{u/4}.
$$
It then follows easily that (for $r$ close to $\bar R$)
$$
u\ge -4\log(\bar R - r) -C.
$$
\qed

\proof{ of Proposition~\ref{p1.7} a)} Given $N\ge 13$, let
$b^{max}$ denote the supremum of all parameters $b\ge-4$ such that
the corresponding extremal solution is singular. { We first
observe that
$$b^{max}>0.$$ In fact, it follows from Sections 5 and 6 that the
extremal solution $u$ associated to parameters $a=b=0$ is strictly
stable :
\begin{align}\label{strictstab}
\inf_{ \varphi \in C_0^\infty(B) } \frac{ \int_{B} (\Delta
\varphi)^2 - \la^* \int_{B} e^{u} \varphi^2}{\int_{B} \varphi^2}
> 0.
\end{align}
Extend $u$ as before on its maximal interval of existence $(0,\bar
R)$. Choosing $R\in (1,\bar R)$ close to $1$, we deduce that
\eqref{strictstab} still holds on the ball $B_R$. In particular,
letting $v(x)=u(Rx)-u(R)$ for $x\in B$, we conclude that $v$ is a
singular stable solution of \eqref{main2} with $a=0$ and
$b=Ru'(R)>0$. By Proposition \ref{240}, we conclude that
$b^{max}>0$.} We now prove that
$$
b^{max}<\infty.
$$
Assume this is not the case and let $u_n$ denote the (singular)
extremal solution associated to $b_n$, where $b_n\nearrow\infty$.
We first observe that there exists $\rho_n\in (0,1)$ such that
$u_n'(\rho_n)=0$. Otherwise, $u_n$ would remain monotone
increasing on $(0,1)$, hence bounded above by $u_n(1)=0$. It would
then follow from \eqref{main2} and elliptic regularity that $u_n$
is bounded. Let $v_n(x)=u_n(\rho_n x)-u_n(\rho_n)$ for $x\in B$
and observe that $v_n$ solves \eqref{main2} with $a=b=0$ and some
$\lambda=\lambda_n$. Clearly $v_n$ is stable and singular. By
Proposition \ref{240}, $v_n$ coincides with $u$, the extremal
solution of \eqref{main2} with $a=b=0$. By standard ODE theory,
$v_n=u$ on $(0,\bar R)$. In addition, $$
b_n=u_n'(1)=\frac1{\rho_n}
v_n'\left(\frac1{\rho_n}\right)=\frac1{\rho_n}
u'\left(\frac1{\rho_n}\right)\to+\infty,
$$ which can only happen
if $ 1/\rho_n\to \bar R$.

Now, since $u_n$ is stable on $B$, $u=v_n$ is stable on
$B_{1/\rho_n}$. Letting $n\to\infty$, we conclude that $u$ is
stable on $B_{\bar R}$. This clearly contradicts Lemma
\ref{lemma71}.

We have just proved that $b^{max}$ is finite. It remains to prove
that $u^*$ is singular when $-4\le b\le b^{max}$. We begin with
the case $b=b^{max}$. Choose a sequence $(b_n)$ converging to
$b^{max}$ and such that the corresponding extremal solution $u_n$
is singular. Using the same notation as above, we find a sequence
$\rho_n\in(0,1)$ such that
$$
\frac1{\rho_n}u'\left(\frac1{\rho_n}\right)=b_n\to b^{max}.
$$
Taking subsequences if necessary and passing to the limit as
$n\to\infty$, we obtain for some $\rho\in(0,1)$
$$
\frac1{\rho}u'\left(\frac1{\rho}\right)=b^{max}.
$$
Furthermore, by construction of $\rho_n$, $u$ is stable in
$B_{1/\rho_n}$ hence in $B_{1/\rho}$. This implies that $v$
defined for $x\in B$ by $v(x)=u(\frac{x}\rho)-u(\frac1\rho)$, is a
stable singular solution of \eqref{main2} with $b=b^{max}$. By
Proposition \ref{240}, we conclude that the extremal solution is
singular when $b=b^{max}$.

When $b=-4$ we have already mentioned in the introduction that
$u^*$ is singular for $N\ge 13$ as a direct consequence of
Proposition~\ref{240} and Rellich's inequality.

So we are left with the case $-4 < b < b^{max}$. Let $u^*_{m}$
denote the extremal solution when $b=b^{max}$, which is singular,
and $\lambda^*_m$ the corresponding parameter. For $0<R<1$ set
\begin{align*}
u_R(x) = u^*_m(R x) - u^*_m(R).
\end{align*}
Then
\begin{align*}
\Delta^2 u_R = \lambda_R e^{u_R} \quad \mbox{where}\quad \lambda_R
= \lambda^*_0 R^4 e^{u^*_m(R)},
\end{align*}
and $u_R = 0$ on $\partial B$, while
\begin{align*}
\frac{d u_R}{d r}(1) = R \frac{d u^*_m}{d r}(R).
\end{align*}
By \eqref{1400}, note that
\begin{align*}
R \frac{d u^*_m}{d r}(R) \to b^{max} \quad \hbox{ as $R\to
1$},\quad\mbox{and}\quad R \frac{d u^*_m}{d r}(R) \to -4 \quad
\hbox{ as $R\to 0$.}
\end{align*}
Thus, for any $-4 < b < b^{max}$ we have found a singular stable
solution to \eqref{main2} (with $a=0$). By Proposition~\ref{240}
the extremal solution to this problem is singular.

%{\color{red} In the proposition it is claimed that $b^{max}>0$. We
%should add a few lines proving this, since Theorem~\ref{t3} only
%guarantees $b^{max} \ge 0$. }

\qed

\medskip

\proof{ of Proposition~\ref{p1.7} b)} Let $b \ge -4$.
Lemma~\ref{bound-u*} applies also for $b\ge -4$ and yields $ u^*
\le \bar u$ where $\bar u(x) = - 4 \log|x|$. We now modify
slightly the proof of Lemma~\ref{L50}. Indeed, consider $w = (4+b)
(1-r^2)/2$ and define
%\begin{align*}
$ u= \bar u - w.$
%\end{align*}
Then
\begin{align*}
\Delta^2 u & = 8 (N-2)(N-4) \frac 1{r^4} = 8 (N-2)(N-4) e^{\bar u}
 = 8 (N-2)(N-4) e^{u + w}
\\
& \le 8 (N-2)(N-4) e^{(4+b)/2} e^u.
\end{align*}
Also $u(1) = 0$, $u'(1) =b$, so $u$ is a subsolution to
\eqref{main2} with parameter $\la_0 = 8 (N-2)(N-4) e^{(4+b)/2}$.

If $N$ is sufficiently large, depending on $b$, we have $\la_0 <
N^2 (N-4)^2 /16$. Then by \eqref{hardy-rellich} $u$ is a stable
subsolution of \eqref{main2} with $\lambda=\lambda_0$. As in
Lemma~\ref{L50} this implies $\lambda^* \le \lambda_0$.

Thus for large enough $N$ we have $ \lambda^* e^{u^*} \le r^{-4}
\, 8 (N-2)(N-4) e^{(4+b)/2} < r^{-4} N^2 (N-4)^2/ 16 $. This and
\eqref{hardy-rellich} show that
\begin{align*}
\inf_{ \varphi \in C_0^\infty(B) } \frac{ \int_{B} (\Delta
\varphi)^2 - \la^* \int_{B} e^{u^*} \varphi^2}{\int_{B} \varphi^2}
> 0
\end{align*}
which is not possible if $u^*$ is bounded. \qed

\medskip

%b) Let $u$ be the extremal solution of \eqref{main2} with boundary
%data $a=b=0$. Then $u$ is singular by Theorem~\ref{t3}. Let $u$
%denote also the unique continuation of this function to a maximal
%interval of existence $0<r<R$ where $R>1$.
%
%We have
%\begin{align}\label{1100}
%u''(1)>0.
%\end{align}
%Indeed, first notice that since $u(1)=u'(1)=0$ and $u(r)>0$ for
%$0<r<1$ we must have $u''(1) \ge 0$. Suppose that $u''(1)=0$.
%Integrating equation \eqref{main2} in $B$ yields
%\begin{align*}
%\int_{\partial B} \pd{\Delta u }{n} = \lambda^* \int_B e^u >0 .
%\end{align*}
%But
%\begin{align*}
%\pd{\Delta u }{n} = u'''(1) + \frac{N-1}{r} u''(1) -
%\frac{N-1}{r^2} u'(1)
%\end{align*}
%and we deduce
%\begin{align*}
%u'''(1)>0.
%\end{align*}
%This would imply $u(r)<0$ for $r<1$ sufficiently close to 1. Thus
%\eqref{1100} holds.
%
%
%\medskip
%
%c) As in the proof of the previous part, let $u$ be the extremal
%solution of \eqref{main2} with zero Dirichlet boundary condition,
%extended to a maximal interval $0<r<R$.
%
%
%
%
%
%We claim that
%\begin{align}\label{1110}
%u'(r)>0 \quad \hbox{for all} 1<r<R.
%\end{align}
%Otherwise define
%\begin{align*}
%r_0 = \inf_{} \{\, 1<r<R \, / \, u'(r)=0 \, \}.
%\end{align*}
%Observe that
%\begin{align*}
%\hbox{$u'(r_0) = 0$, $u'(r)>0$ and $u(r)>0$ for $1<r<r_0$.}
%\end{align*}
%Then the function
%\begin{align*}
%v(r) = u(r r_0) - u (r_0)
%\end{align*}
%is a solution to \eqref{main2} with parameter $\lambda =
%\lambda^*_0 r_0^4 e^{u(r_0)}$ and zero Dirichlet data. By the
%maximum principle $v>0$ in $B$, which is not possible, since
%$u(r)<u(r_0)$ for $1<r<r_0$. This establishes \eqref{1110}.

\bigskip

\noindent{\bf Acknowledgements.}

\medskip
This work has been partly supported by grants FONDECYT 1050725 and
FONDAP, Chile. Part of it was done when J.D. was visiting
Laboratoire Ami\'enois de Math\'ematique Fondamentale et
Appliqu\'ee at Facult\'e de Math\'ematique et
d'In\-for\-ma\-tique, Amiens, France and when L.D. was visiting
Departamento de Ingenier\'{\i}a Matem\'atica at Universidad de
Chile, Santiago, Chile. They are grateful for the hospitality and
support. I. Guerra was supported by FONDECYT 1061110. M.
Montenegro was supported by Project CNPq-CNRS 490093/2004-3.

\medskip
J.D. and L.D. are indebted to J. Coville for useful conversations
on some of the topics in the paper.

\medskip
The authors thank the referee for his careful reading of the
manuscript.

\end{document}